\documentclass[%
 aip,
 amsmath,amssymb,
 reprint,%
]{revtex4-1}

\usepackage{graphicx}
\usepackage{dcolumn}
\usepackage{bm}

\usepackage[utf8]{inputenc}
\usepackage[T1]{fontenc}
\usepackage{mathptmx}
\usepackage{etoolbox}
\usepackage{url}

\usepackage{color}

\usepackage{graphicx}
\usepackage{subfigure}
\usepackage{overpic} 
\usepackage{stfloats}
\usepackage{hyperref} 

\makeatletter
\def\@email#1#2{%
 \endgroup
 \patchcmd{\titleblock@produce}
  {\frontmatter@RRAPformat}
  {\frontmatter@RRAPformat{\produce@RRAP{*#1\href{mailto:#2}{#2}}}\frontmatter@RRAPformat}
  {}{}
}%
\makeatother
\begin{document}

\preprint{AIP/123-QED}

\title[Fractional SEIR Model and Data-Driven Predictions of COVID-19 Dynamics of Omicron Variant]{Fractional SEIR Model and Data-Driven Predictions of COVID-19 Dynamics of Omicron Variant}
\author{Min Cai}
\affiliation{Department of Mathematics, Shanghai University, 99 Shangda Road, Shanghai 200444, China} 
\author{George Em Karniadakis}
\affiliation{Division of Applied Mathematics, Brown University, 170 Hope Street, Providence,
RI 02906, USA}
\author{Changpin Li}
\affiliation{Department of Mathematics, Shanghai University, 99 Shangda Road, Shanghai 200444, China}   
\email{lcp@shu.edu.cn}

\date{\today}

\begin{abstract}
We study the dynamic evolution of COVID-19 cased by the Omicron variant via a fractional susceptible-exposed-infected-removed (SEIR) model.  Preliminary data suggest that the symptoms of Omicron infection are not prominent and the transmission is therefore more concealed, which causes a relatively slow increase in the detected cases of the new infected at the beginning of the pandemic. To characterize the specific dynamics, the Caputo-Hadamard fractional derivative is adopted to refined the classical SEIR model. Based on the reported data, we infer the fractional order, time-dependent parameters, as well as unobserved dynamics of the fractional SEIR model via fractional physics-informed neural networks (fPINNs). Then, we make short-time predictions using the learned fractional SEIR model. 
\end{abstract}

\maketitle

\begin{quotation}
The susceptible-exposed-infected-removed (SEIR) model is an epidemiological model that characterizes transmission of infectious diseases with incubation. Introducing fractional differential operators to modify the model provides an alternative approach as the history-dependence of the epidemic can be taken into account. Since the transmission of the Omicron variant is more concealed than earlier variants of COVID-19 including Beta variant and Delta variant, early reported data of new infected cases reveal a relatively slow increase instead of fast increase, the Caputo-Hadamard derivative is employed in the fractional SEIR model. To calibrate the model, fractional physics-informed neural networks (fPINNs) which is a novel deep learning framework solving data-driven forward and inverse problems of fractional differential equations, is applied. Then, the calibrated fractional SEIR model is able to give reliable short-term predictions of the COVID-19 dynamics caused by Omicron variant to a certain extent.  
\end{quotation}

\section{\label{sec:Introduction}Introduction} 

Since the end of 2021, outbreaks of the epidemic caused by the Omicron variant have  occurred in numerous cities around the world, such as Berlin, New York, and Shanghai. This new variant spreads more easily than earlier variants (such as Beta variant and Delta variant) of the virus that cause COVID-19. In addition, it was reported that anyone with the Omicron infection, regardless of vaccination status or whether or not they have symptoms, can spread the virus to others \cite{CDC}. The outbreak caused by Omicron variant spreads rapidly and surpassed the outbreak of the COVID-19 at the beginning of 2020. All sectors of the society are taking great efforts to fight this epidemic and measures in many aspects keep  fine-tuning during the epidemic. Therefore, a whole picture of the transmission of this disease infection and reliable predictions of this epidemic are of great significance. 

Epidemiological models offer computationally expedient predictions for quantifying the dynamics of an epidemic. Although several compartmental models were developed to study the epidemic of COVID-19 occurred at the beginning of 2020, it is necessary to construct a modified model describing the transmission of the disease infection caused by the new variant. In this paper, we modify the susceptible-exposed-infected-removed (SEIR) epidemiological model by introducing the Caputo-Hadamard fractional derivative to characterize the epidemic of COVID-19 caused by the Omicron variant, mainly due to the following reasons: 
\begin{itemize}
    \item The Omicron variant has an incubation period and SEIR model is a basic model suitable for the infectious disease with incubation; 
    \item The reported data usually include the daily new infected cases and the daily new removed cases (recovered cases and dead cases); hence, an epidemiological model consisting of $S$, $E$, $I$ and $R$ compartments is necessary;  
    \item As the transmission of Omicron variant is more concealed, the detected cases of Omicron infections have a relatively slow increase at the beginning of the outbreak. The Caputo-Hadamard derivative whose integration kernel is given by the power of the logarithmic function is therefore suitable to characterize this specific process. 
\end{itemize}

To draw a complete picture of the fractional SEIR model, the model parameters should be known. Methods of calibrating those parameters via reported data include the particle swarm optimization (PSO) algorithm \cite{He2020}, nonlinear least square curve fitting approach \cite{Oud2021}, and machine learning \cite{EC2021}.  The former two approaches are usually applied to calibrating the time-independent parameters while the third can be applied to inferring both time-dependent and time-independent variables. In this paper, a deep learning framework called fractional physics-informed neural networks (fPINNs) is applied to solving the inverse problem corresponding to the fractional SEIR model. We obtain the fractional order, the time-dependent model parameters, and the compartments without reported data. Once the fractional model is learned, this specific model can be employed to predict the epidemic by  solving the forward problem. Though the measures to control the epidemic keep adjusting, which may result in uncertainty, simulations illustrate that the fPINNs is well-trained and the learned fractional SEIR model is able to forecast the short-term dynamics to a certain extent. 

The paper is organized as follows: In Sec. \ref{sec:Methodology} we introduce the fractional SEIR model and the corresponding fPINNs as well as its implementation details. 
In Sec. \ref{sec:Results} we present simulation results based on the reported data from Shanghai  and in Sec. \ref{sec:Discussion} we include a discussion and concluding remarks.

\section{\label{sec:Methodology}Methodology}
\subsection{\label{sec:fSEIR}Fractional SEIR Model}
In compartmental models, the population is divided into compartments with labels, such as $S$, $E$, $I$, and $R$ corresponding to Susceptible, Exposed, Infectious, and Removed compartments. People may progress between the  compartments under certain transmission principles. The following SEIR model is suitable for the infectious disease with incubation:
\begin{equation*}
    \begin{aligned}
    &\frac{\rm d}{{\rm d}t}S(t)
    = -\frac{\beta(t)S(t)I(t)}{N}, 
    \\
    &\frac{\rm d}{{\rm d}t}E(t)
    =\frac{\beta(t)S(t)I(t)}{N} - \sigma E(t),
    \\
    &\frac{\rm d}{{\rm d}t}I(t)
    =\sigma E(t) - \mu(t)I(t),
    \\
    &\frac{\rm d}{{\rm d}t}R(t)
    =\mu(t)I(t).
    \end{aligned}
\end{equation*}
Here $N$, $S(t)$, $E(t)$, $I(t)$, and $R(t)$ represent the population of the whole society, susceptible people, exposed people, infected people, and removed (including recovered and dead) people, respectively. The model parameters $\beta(t)$, $\mu(t)$, and $\sigma$ denote the rate of transmission for the susceptible to exposed, the rate of transmission for the exposed to infected, and the incubation rate. Here $\sigma^{-1}$ is the value of incubation period \cite{Yang2020}. Since the average incubation period of Omicron is about 3 days \cite{CDC}, we assumed that $\sigma = \frac{1}{3}$. 

There are several variants of the SEIR model to refine the model. Yang et al. \cite{Yang2020} propose a modified SEIR model by introducing the move-in and move-out parameters to the transmission between the  susceptible people and the exposed people. Some other existing literature modified the SEIR model by adding more compartments into the original model \cite{Hao2020,Zhang2021}. 

Herein, we modify the SEIR model for the epidemic of COVID-19 caused by Omicron variant via introducing the fractional time derivative. Fractional calculus can be very expressive across different fields, as the fractional integral equation and fractional differential equation are able to characterize materials and phenomena with nonlocality and history dependency\cite{LiCai2019}. Therefore, compartmental models cooperating with fractional differential operators offer a novel option to model infection diseases \cite{EC2021,Tuan2020,Xu2020}. In most of the existing works the Caputo derivative is used to modify the models due to the fast increase of the reported data on daily new infected as was the case at the beginning of the epidemic in 2020. However, the detected cases of Omicron infections have a relatively slow increase at the beginning of the outbreak, mainly due to the more concealed  transmission of the new variant. The Caputo-Hadamard fractional derivative is more suitable in describing the ultra slow process \cite{Denisov2010,LiLi2021,Lomnitz1956,Lomnitz1957,Lomnitz1962,Mainardi2017}, mainly because its integration kernel is given by the power of logarithmic function instead of the power function. Therefore, we adopt the Caputo-Hadamard derivative here to model the epidemic of COVID-19 caused by Omicron variant, and the resulting fractional SEIR model is given by: \begin{equation*}
    \begin{aligned}
    &{}_{CH}{\rm D}^{\alpha}_{1,t}S(t)
    = -\frac{\beta S(t)I(t)}{N}, 
    \\
    &{}_{CH}{\rm D}^{\alpha}_{1,t}E(t)
    =\frac{\beta S(t)I(t)}{N} - \sigma E(t)
    \\
    &{}_{CH}{\rm D}^{\alpha}_{1,t}I(t)
    =\sigma E(t) - \mu I(t)
    \\
    &{}_{CH}{\rm D}^{\alpha}_{1,t}R(t)
    =\mu I(t)
    \\
    &{}_{CH}{\rm D}^{\alpha}_{1,t}I^{c}(t)
    =\sigma E(t). 
    \end{aligned}
\end{equation*}
The last equation in the above system follows from the fact that the cumulative infected cases $I^{c}(t)$ equals to the summation of the inflows of the $I(t)$ compartment. The daily new infected cases $I^{n}(t)$ at $t= t_{j}$ expresses the difference $I^{c}(t_{j+1})-I^{c}(t_j)$. When the data for daily new infected cases are available, the equation containing the variable $I^{c}(t)$ is needed for training the neural networks. Similarly, as there is no outflow for the $R(t)$ compartment, the daily new removed cases can be computed as $R^{n}(t_{j}) = R(t_{j+1})-R(t_j)$. Here the Caputo-Hadamard derivative ${}_{CH}{\rm D}^{\alpha}_{a,t}$ is defined as 
\begin{equation*}
     {}_{CH}{\rm D}^{\alpha}_{a,t}u(t)
    =\frac{1}{\Gamma(n-\alpha)}\int_{a}^{t}
    \left(\log\frac{t}{s}\right)^{n-\alpha-1}\delta^{n}u(s)
    \frac{{\rm d}s}{s},
    \ 
    a>0, 
\end{equation*}
where $\delta^{n}u(s)=\left(s\frac{\rm d}{{\rm d}s}\right)^{n}u(s)
=\delta(\delta^{n-1}u(s))$, $\delta^{0}u(s)=u(s)$, $n-1<\alpha<n\in\mathbb{Z}^{+}$.

\subsection{\label{sec:fPINNs}Fractional Physics-informed Neural Networks and the Implementation}
 
Physics-informed neural networks (PINNs) \cite{RP2019} and fractional physics-informed neural networks (fPINNs) \cite{PL2019} were introduced by Karniadakis and his collaborators. These two deep learning frameworks have been successfully applied in solving forward and inverse problems in many practical applications that can be characterized by integer-order differential equations and fractional differential equations \cite{JK2020,MJ2020,RB2019,YZ2020}. 

The basic idea of PINNs and fPINNs is training the neural networks to solve semi-supervised learning tasks while respecting any given laws of physics described by general differential equations. That is to say, the trained neural networks should not only fit the training data well, but also satisfy the governing differential equations. In the fPINNs formulation of the present paper, separate deep neural networks with input $t$ are adopted to represent the 
compartments $S(t)$, $E(t)$, $I(t)$, $R(t)$, $I^{c}(t)$ as well as the time-dependent parameters $\beta(t)$, $\mu(t)$. Each network is parameterized by a set of parameters $\Theta$ as weights and biases of the network. The resulting deep neural networks are denoted by $S_{NN}(t)$, $E_{NN}(t)$, $I_{NN}(t)$, $R_{NN}(t)$, $I^{c}_{NN}(t)$, $\beta_{NN}(t)$, and $\mu_{NN}(t)$. The loss function of fPINNs for the considered fractional SEIR model is defined as   
\begin{equation*}
\begin{aligned}
      L(\Theta)
    =&{\rm MSE}_{u} + {\rm MSE}_{r}, 
\end{aligned}
\end{equation*}
with ${\rm MSE}_{u}$ and ${\rm MSE}_{r}$ being the losses from the data fitting and the fractional ordinary differential equation (ODE) residual. 

Most of the data sets consist of daily new infected cases and daily new removed (recovered and dead) cases while there is no information on the amount of susceptible people $S(t)$ and the amount of exposed people $E(t)$. Let us denote by $\{t_{u}^{j}\}_{j = 1}^{N_u}$ the temporal nodes for the training data with $N_{u}\in\mathbb{Z}^{+}$ being the number of training data. In the present paper, $\{t_{u}^{j}\}_{j = 1}^{N_u} = \{1, 2, \ldots, N_u\}$. It can be derived that $I^{c}(t_{u}^{j+1}) = I^{c}(t_{u}^{j}) + I^{n}(t_{u}^{j+1})$, $R(t_{u}^{j+1}) = R(t_{u}^{j}) + R^{n}(t_{u}^{j+1})$, and $I(t_{u}^{j}) = I^{c}(t_{u}^{j}) - R(t_{u}^{j})$. As a result, the loss function for data fitting is given by: 
\begin{equation*}
    \begin{aligned}
           {\rm MSE}_{u}
         &=\sum\limits_{j=1}^{N_{u}}\bigg|I_{NN}^{n}(t_{u}^{j})-I^{n}(t_{u}^{j})\bigg|^{2} 
         +\sum\limits_{j=1}^{N_{u}}\bigg|I_{NN}^{c}(t_{u}^{j})-I^{c}(t_{u}^{j})\bigg|^{2}
         \\
         &+\sum\limits_{j=1}^{N_{u}}\bigg|R_{NN}^{n}(t_{u}^{j})-R^{n}(t_{u}^{j})\bigg|^{2} 
         +\sum\limits_{j=1}^{N_{u}}\bigg|R_{NN}(t_{u}^{j})-R(t_{u}^{j})\bigg|^{2}
         \\
         &+\sum\limits_{j=1}^{N_{u}}\bigg|I_{NN}(t_{u}^{j})-I(t_{u}^{j})\bigg|^{2}. 
    \end{aligned}
\end{equation*}

For the considered fractional SEIR model, the fractional ODE residual can be defined by 
\begin{equation*}
    \begin{aligned}
           {\rm MSE}_{r}
         &=\sum\limits_{j=1}^{N_{r}}\bigg|{}_{CH}{\rm D}^{\alpha}_{1,t}S_{NN}(t_{r}^{j})+\frac{\beta S_{NN}(t_{r}^{j})_{NN}I(t_{r}^{j})}{N}\bigg|^{2}
         \\
         &+\sum\limits_{j=1}^{N_{r}}\bigg|{}_{CH}{\rm D}^{\alpha}_{1,t}E_{NN}(t_{r}^{j})
         -\frac{\beta S_{NN}(t_{r}^{j})I_{NN}(t_{r}^{j})}{N}+\sigma E_{NN}(t_{r}^{j})\bigg|^{2}
         \\
         &+\sum\limits_{j=1}^{N_{r}}\bigg|{}_{CH}{\rm D}^{\alpha}_{1,t}I_{NN}(t_{r}^{j})
         -\sigma E_{NN}(t_{r}^{j}) + \mu I_{NN}(t_{r}^{j})\bigg|^{2}
         \\
         &+\sum\limits_{j=1}^{N_{r}}\bigg|{}_{CH}{\rm D}^{\alpha}_{1,t}R_{NN}(t_{r}^{j})
         -\mu I_{NN}(t_{r}^{j})\bigg|^{2}
         \\
         &+\sum\limits_{j=1}^{N_{r}}\bigg|{}_{CH}{\rm D}^{\alpha}_{1,t}I^{c}_{NN}(t_{r}^{j})
         -\sigma E_{NN}(t_{r}^{j})\bigg|^{2}
    \end{aligned}
\end{equation*}
with $N_{r}\in\mathbb{Z}^{+}$ being the number of the residual points. 

It should be noted that while calculating the fractional ODE residual ${\rm MSE}_{r}$, Caputo-Hadamard derivatives of neural networks approximating the $S(t)$, $E(t)$, $I(t)$, $R(t)$, and $I^{c}(t)$ compartments should be evaluated. Different from the integer-order case in which the derivatives of the neural networks with respect to the input $t$ and  parameters $\Theta$ can be computed by applying the chain rule for differentiating compositions of functions and can be realized by the automatic differentiation \cite{baydin2017automatic}, fractional derivatives of the neural networks should be obtained via numerical approximation as the chain rule no longer holds for the fractional differentiation. Numerical methods for fractional ordinary differential equations with Caputo-Hadamard derivative can be found in the work by Gohar et al. \cite{GoharLiLi2020,GoharLiYin2020}. In this paper, the following numerical formula proposed by Fan et al. \cite{FanLi2022} is adopted to evaluate the Caputo-Hadamard derivative,   
\begin{equation*}
    \begin{aligned}
     &\left.{}_{CH}{\rm D}^{\alpha}_{a,t}u(t)\right|_{t=t_{k}} 
    \approx
    \frac{1}{\Gamma(2-\alpha)}
    \sum\limits_{j=1}^{k}c_{j,k}^{(\alpha)}(u^{j}-u^{j-1})
    \\
    =
    &\frac{1}{\Gamma(2-\alpha)}
    \left(
    c_{k,k}^{(\alpha)}u^{k}
    +\sum\limits_{j=1}^{k-1}(c_{j,k}^{(\alpha)}-c_{j+1,k}^{(\alpha)})u^{j}
    -c_{1,k}^{(\alpha)}u^{0}
    \right),
    \end{aligned}
\end{equation*}
where $a\leq t_0<\cdots<t_{k}<\cdots<t_{N}$ and 
\begin{equation*}
     c_{j,k}^{(\alpha)}
    =\left\{
    \begin{array}{ll}
    \frac{1}{\log\frac{t_1}{t_0}}\left(a_{1,k}^{(\alpha)}+b_{2,k}^{(\alpha)}\right), 
    & j=1, 
    \\
    \frac{1}{\log\frac{t_j}{t_{j-1}}}
    \left(a_{j,k}^{(\alpha)}-b_{j,k}^{(\alpha)}+b_{j+1,k}^{(\alpha)}\right), 
    & 2\leq j\leq k-1, 
    \\
    \frac{1}{\log\frac{t_k}{t_{k-1}}}\left(a_{k,k}^{(\alpha)}-b_{k,k}^{(\alpha)}\right), 
    & j=k,  
    \end{array}
     \right.
\end{equation*}
with 
\begin{equation*}
    \begin{aligned}
     a_{j,k}^{(\alpha)}
    =\left(\log\frac{t_k}{t_{j-1}}\right)^{1-\alpha}
    -\left(\log\frac{t_k}{t_j}\right)^{1-\alpha},
    \end{aligned}
\end{equation*}
and  
\begin{equation*}
\begin{aligned}
     b_{j,k}^{(\alpha)}
    =&\frac{1}{\log\frac{t_j}{t_{j-2}}}
    \left\{\log\frac{t_j}{t_{j-1}}\left[\left(\log\frac{t_k}{t_j}\right)^{1-\alpha}+\left(\log\frac{t_k}{t_{j-1}}\right)^{1-\alpha}\right]\right.
    \\
     &\left.+\frac{2}{2-\alpha}\left[\left(\log\frac{t_k}{t_j}\right)^{2-\alpha}-\left(\log\frac{t_k}{t_{j-1}}\right)^{2-\alpha}\right]\right\}.
\end{aligned}
\end{equation*}   

\section{\label{sec:Results}Results}
We are interested to infer the unobserved dynamics $S(t)$, $E(t)$ and the parameters $\alpha$, $\beta(t)$, $\mu(t)$ by solving the inverse problem of the fractional SEIR model with Caputo-Hadamard derivatives. Subsequently, given the learned model, we will attempt to forecast the epidemic dynamics by solving the forward problem. The simulation in this section is based on the Shanghai data set as the Shanghai data set include exhaustive information on the $I(t)$ compartment and the $R(t)$ compartment. Simulations based on other data sets can be implemented similarly and may infer different model parameters for different areas. 

\subsection{Basic Settings for Simulation}
The basic settings for the implementation of fPINNs in the simulation are as follows:
\begin{itemize}
    \item Reported data: The data is collected from the published data from the National Health Commission of the People's Republic of China \cite{EpidemicNotification} since 27 Feb 2022. The original reported data include the daily new infected cases, the daily new recovered cases, and the daily new dead cases, from which the cumulative infected cases, the current infected cases, the daily new removed cases (sum of the daily new recovered cases and the daily new dead cases), and the cumulative removed cases (also the current removed cases) can be obtained via direct calculation. 
    \item Preconditioned data: Taking the effect of the hysteresis on the reported data, the training data is given by the seven-days average of the reported data. 
    \item Structure of the neural network: Neural networks with 5 hidden layers and 20 neurons per layer are employed to approximate the compartments $S(t)$, $E(t)$, $I(t)$, $R(t)$, and $I^{c}(t)$. The neural networks with a single hidden layer and 5 neurons per layer are employed to approximate the time-dependent model parameters $\beta(t)$ and $\mu(t)$.
    \item Residual points: In the formulation of PINNs corresponding to the integer-order differential equation, the computation of derivatives for the neural networks can be realized by applying the automatic differentiation technique and the derivatives at arbitrary point can be obtained through this technique. Therefore, the residual points in the PINNs formulation can be chosen arbitrarily. However, residual points in fPINNs formulation are usually grid points in the discretization of the fractional derivative as  fractional derivatives of neural networks are obtained by numerical approaches. The partition of the computational domain affects the approximation of the fractional derivative as the numerical error is a function of the stepsize. In the aforementioned numerical formula for Caputo-Hadamard derivative, the numerical error is of $\mathcal{O}(\tau^{3-\alpha})$ for the uniform mesh $t_{k} = 1+k\tau \in [1, N_u]$ with $\tau = \frac{N_{u}-1}{N_{r}}$ being the time step. In the simulation of this section, we set $\tau = \frac{1}{10}$ and the residual points are therefore chosen as $t_{r}^{j} = 1+\frac{(j-1)(N_{u}-1)}{N_{r}} = 1 + \frac{j-1}{10}, j = 1,2,\ldots, N_r$.  
\end{itemize}

\subsection{\label{sec:Fitting} Data Fitting}
\begin{figure*}[h] 
\centering
\subfigure[\label{Fig:New_I}New infected cases $I^{n}(t)$]{\includegraphics[width = 8cm]{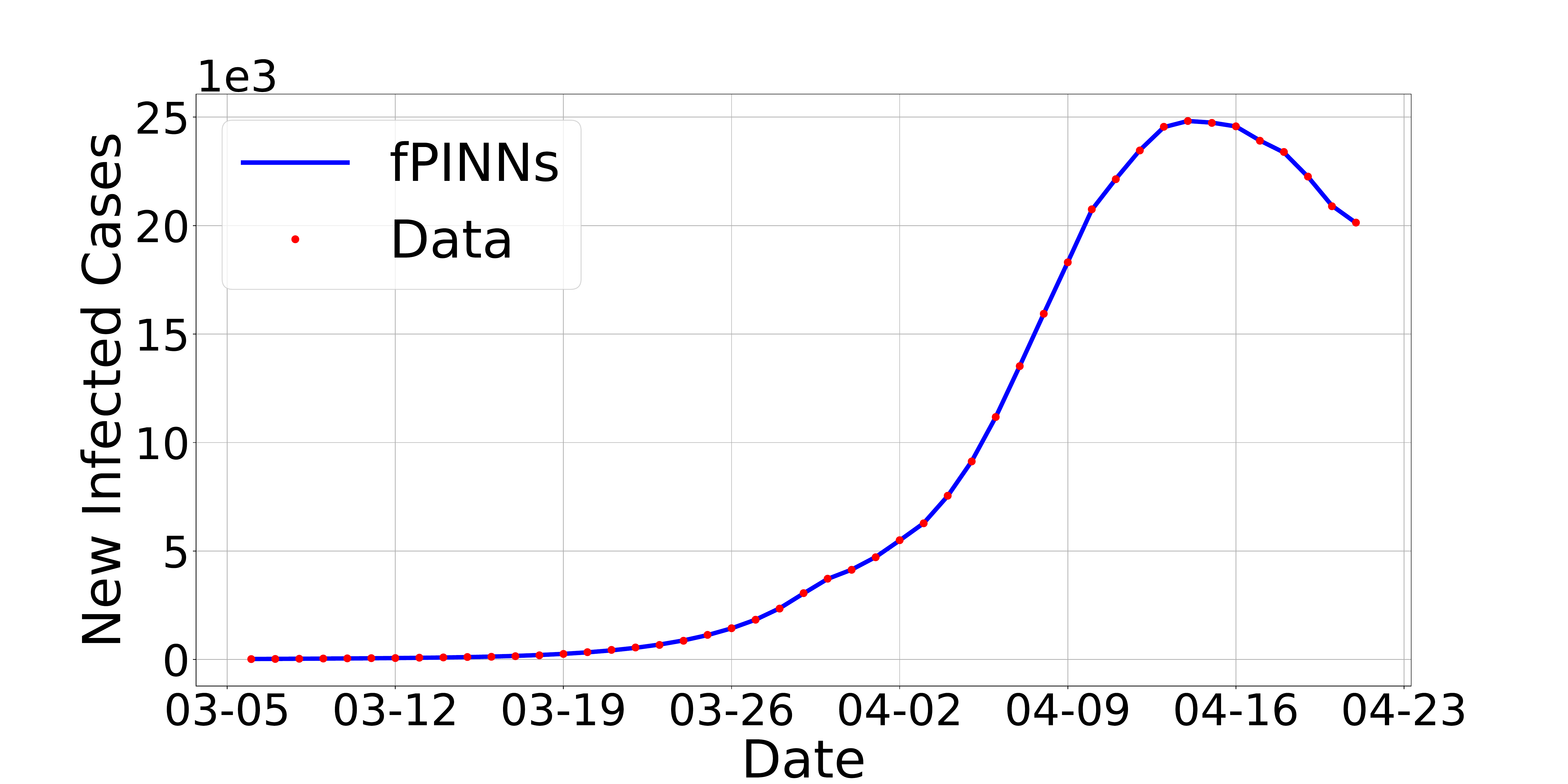}}
\subfigure[\label{Fig:Cumulative_I}Cumulative infected cases $I^{c}(t)$]{\includegraphics[width = 8cm]{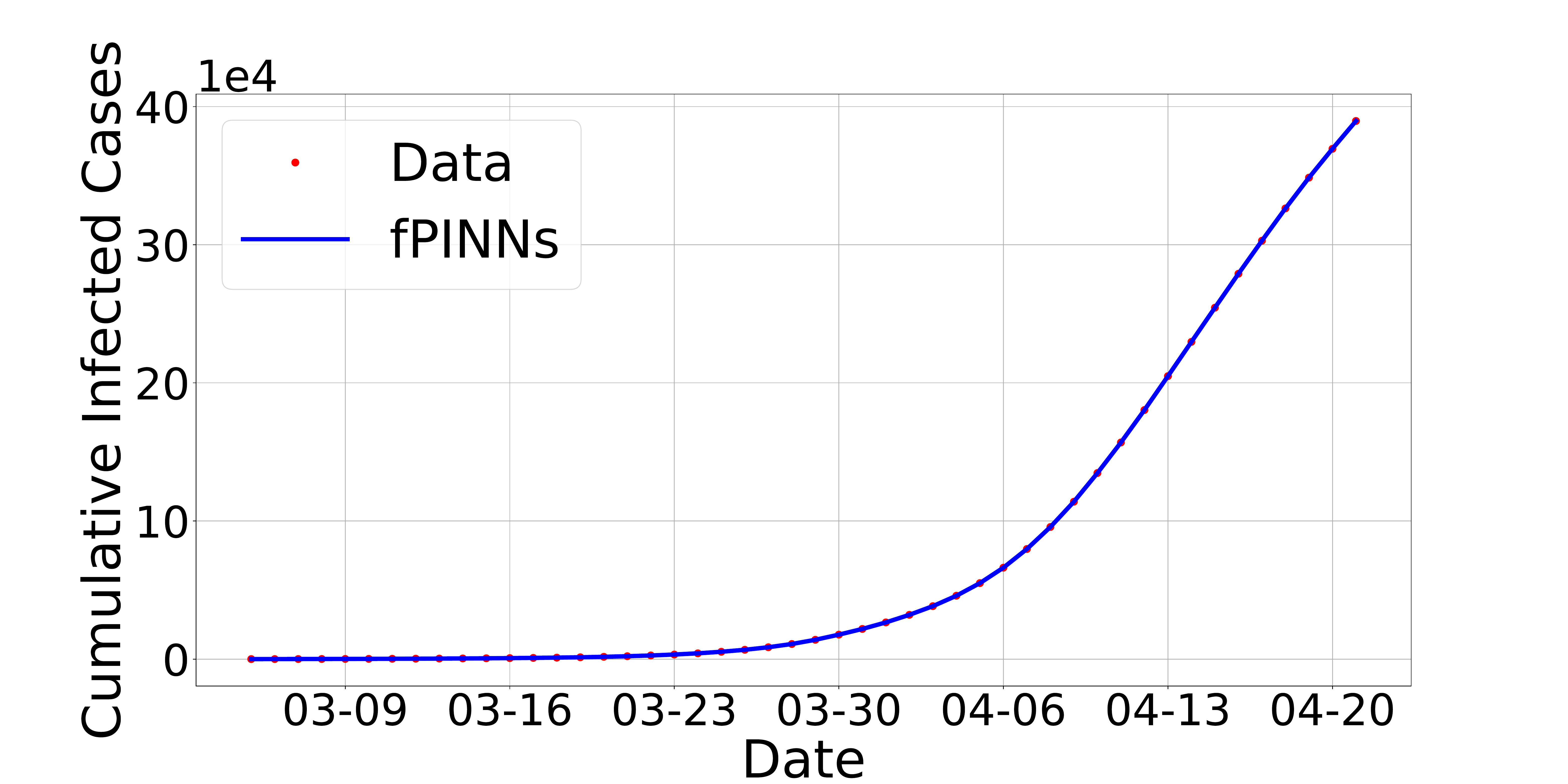}}
\\
\centering
\subfigure[\label{Fig:New_R}New removed cases $R^{n}(t)$]{\includegraphics[width = 8cm]{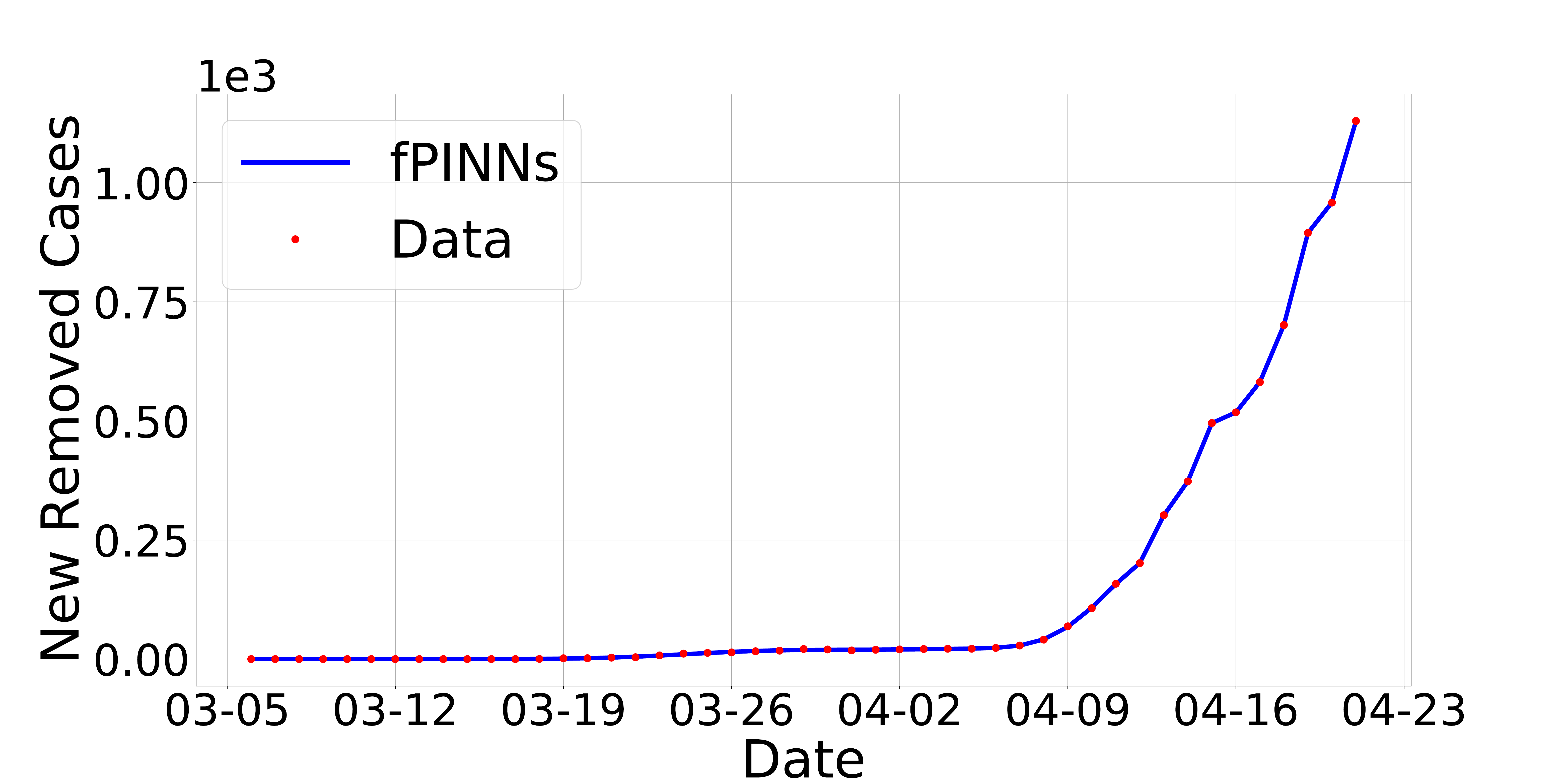}}
\subfigure[\label{Fig:Current_R}Current removed cases $R(t)$]{\includegraphics[width = 8cm]{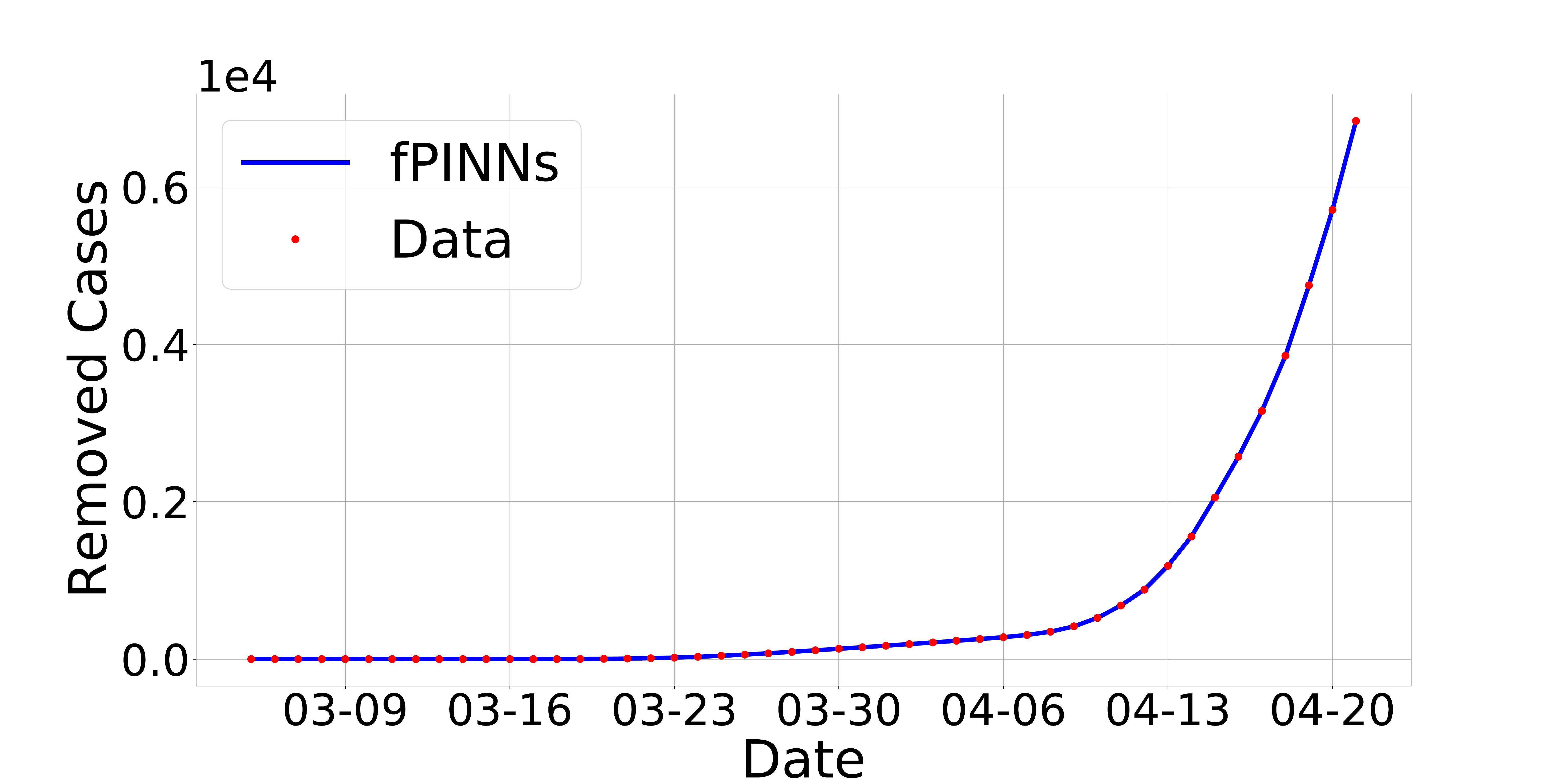}}
\\
\centering
\subfigure[\label{Fig:Current_I}Current infected cases $I(t)$]{\includegraphics[width = 8cm]{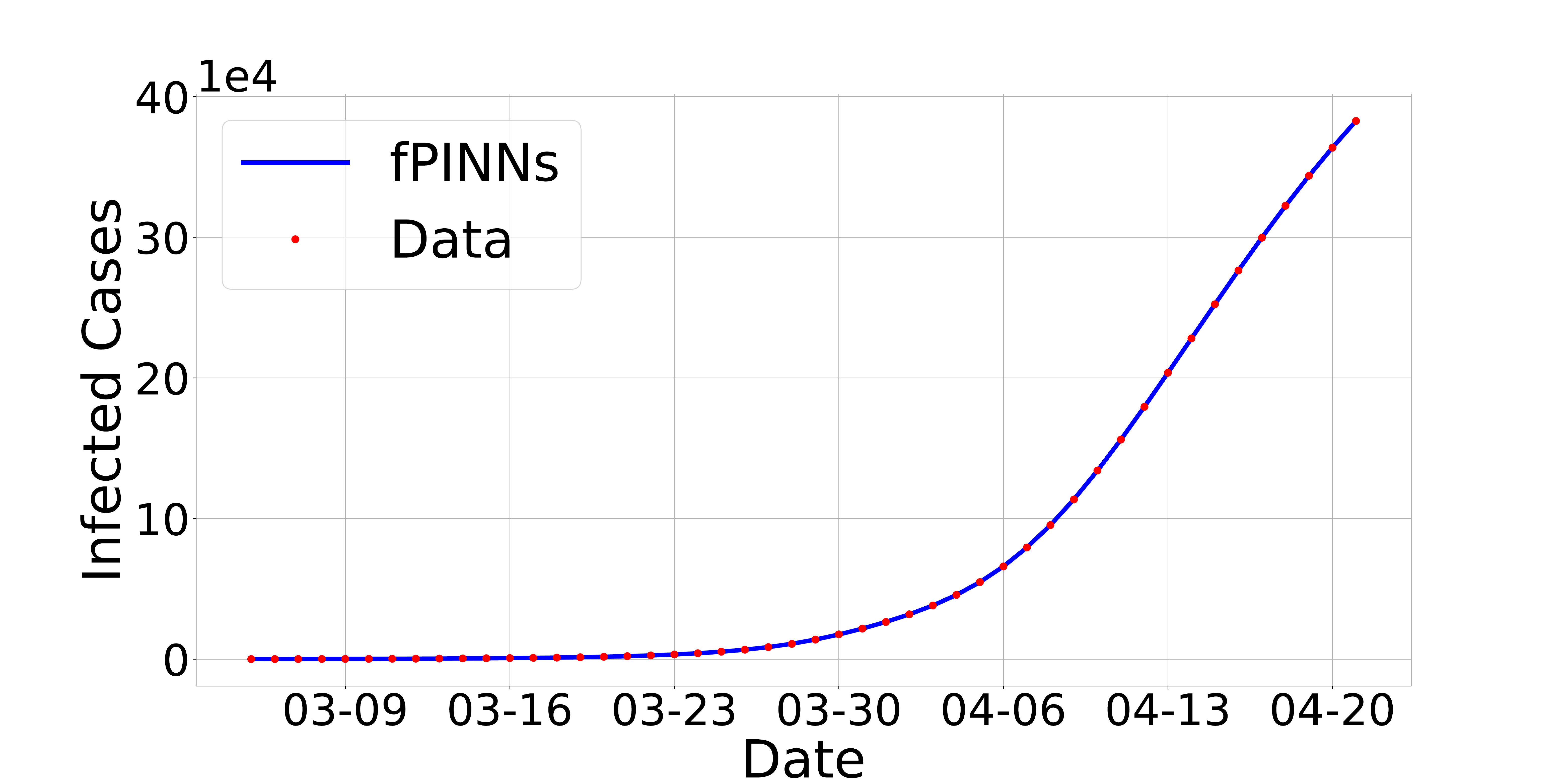}}
\caption{\label{Fig:DataFitting}fPINNs results: Data fitting based on the reported data from 27 February 2022 to 21 April 2022.}
\end{figure*}

In this part, fPINNs results on data fitting for the new infected cases $I^{n}(t)$, the cumulative infected cases $I^{c}(t)$, the current infected cases $I(t)$, the new removed cases $R^{n}(t)$, and the current removed cases $R(t)$ are presented based on the reported data from 27 February 2022 to 21 April 2022. We precondition the reported data by taking seven-days average. The corresponding training data is therefore from 5 March 2022 to 21 April 2022. We observe in Fig. \ref{Fig:DataFitting} that the fPINNs inferences coincide with the ground truth, which indicates that the neural networks are well-trained. 

\subsection{\label{sec:Inferences}Inferences}

\begin{figure*}[h] 
\centering
\subfigure[\label{Fig:Current_S}Current susceptible cases $S(t)$]{\includegraphics[width = 8cm]{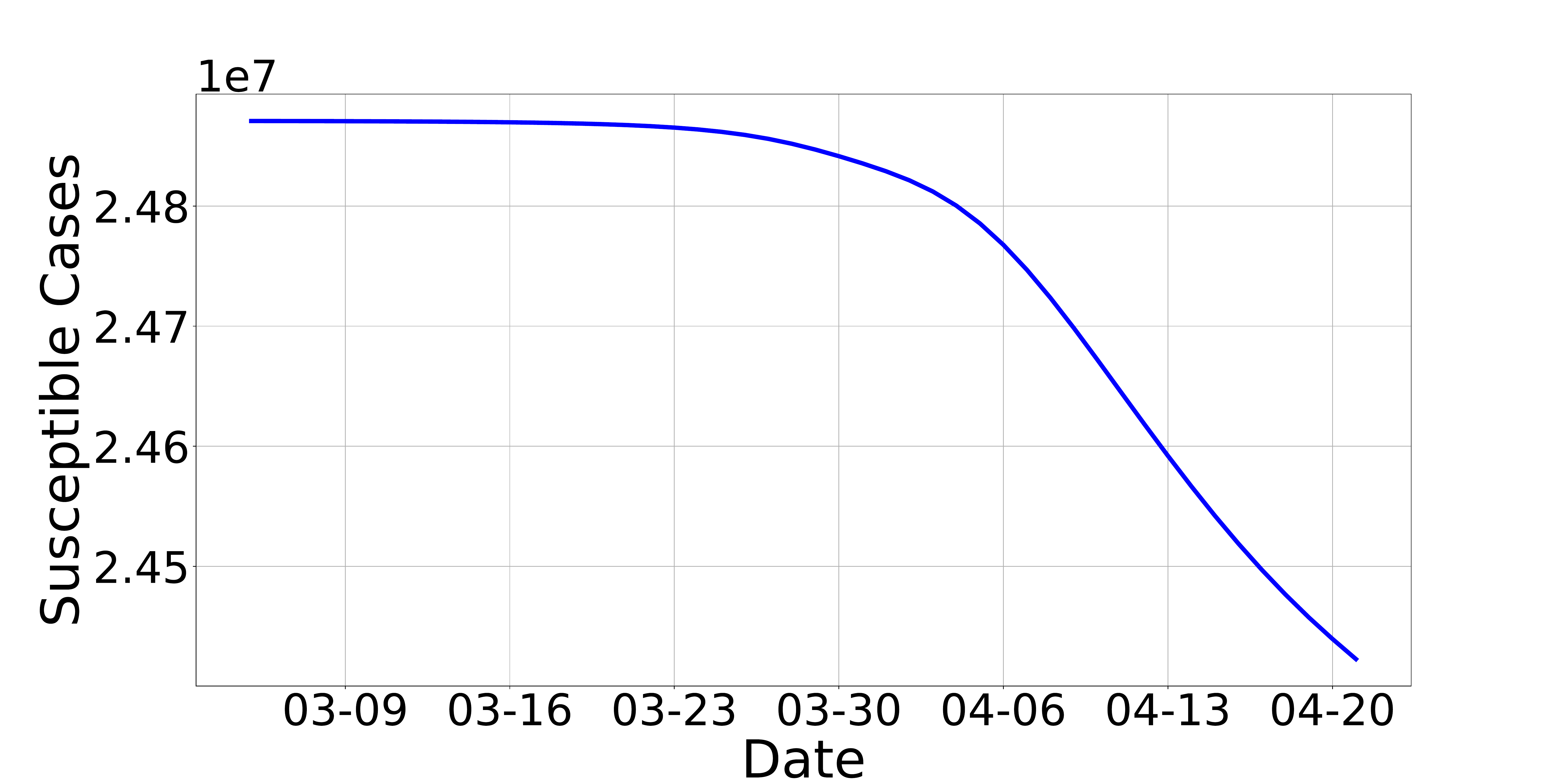}}
\subfigure[\label{Fig:Current_E}Current exposed cases $E(t)$]{\includegraphics[width = 8cm]{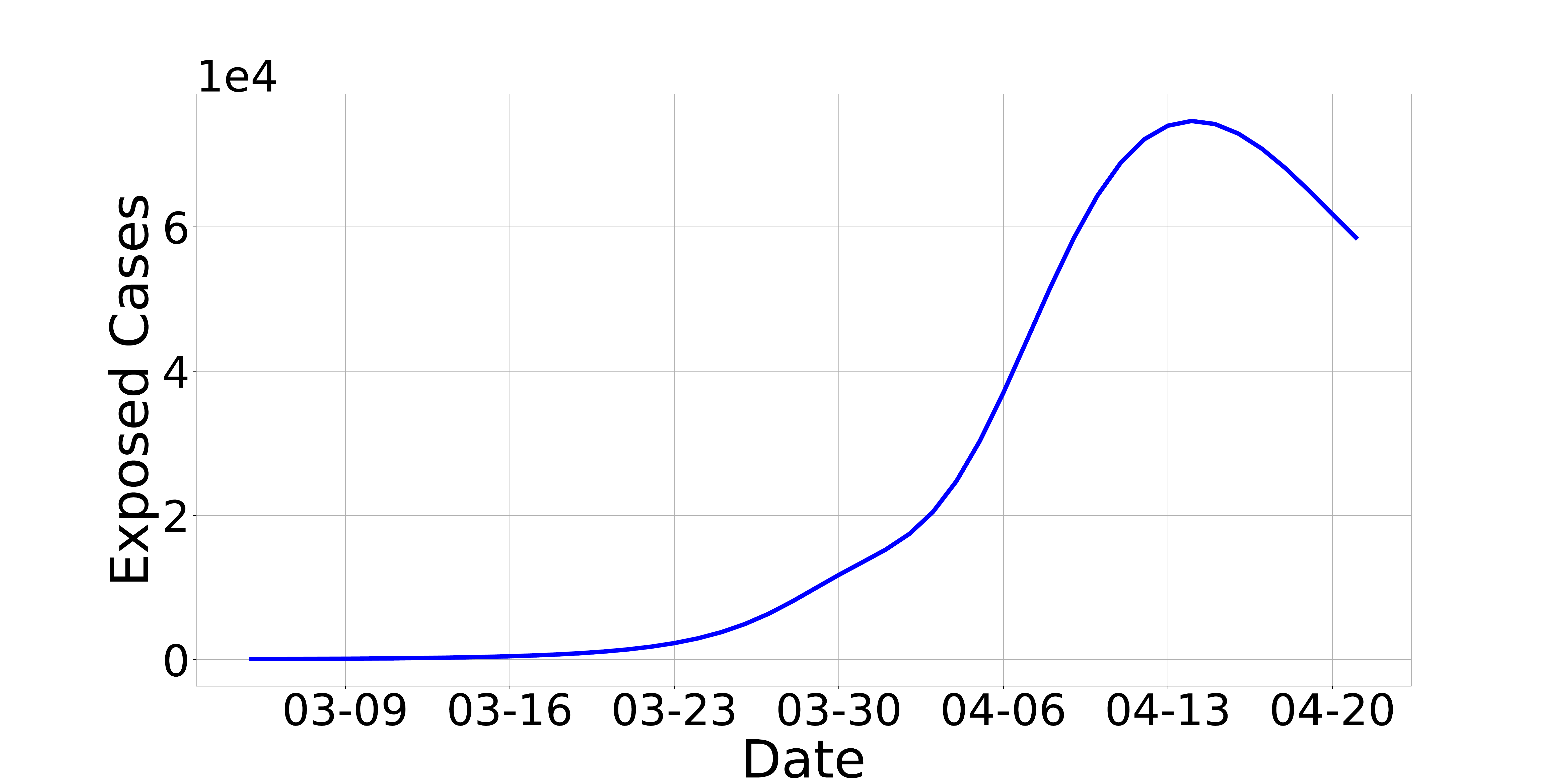}}
\\ 
\centering
\subfigure[\label{Fig:Beta}Model parameter $\beta(t)$]{\includegraphics[width = 8cm]{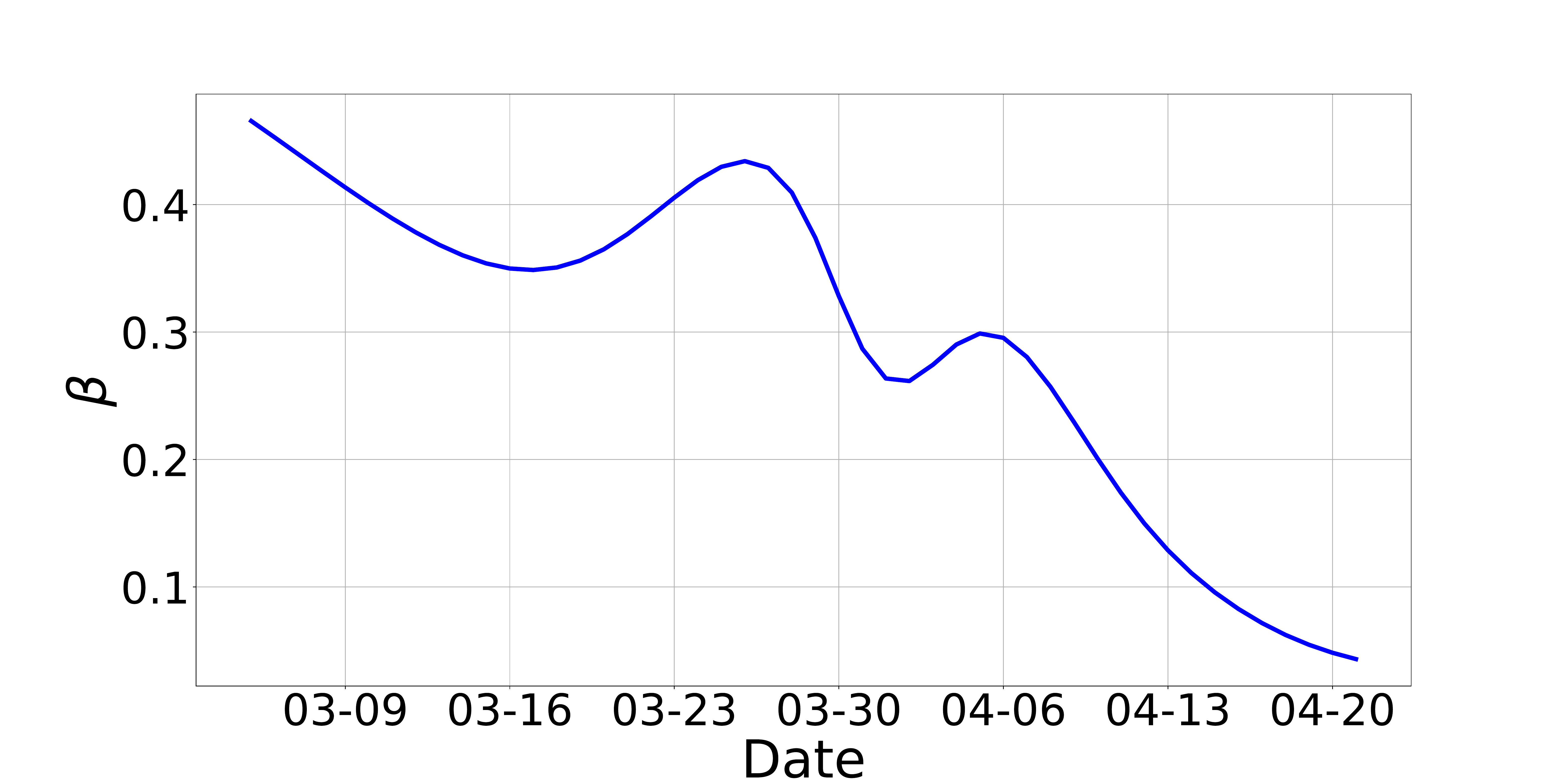}}
\subfigure[\label{Fig:mu}Model parameter $\mu(t)$]{\includegraphics[width = 8cm]{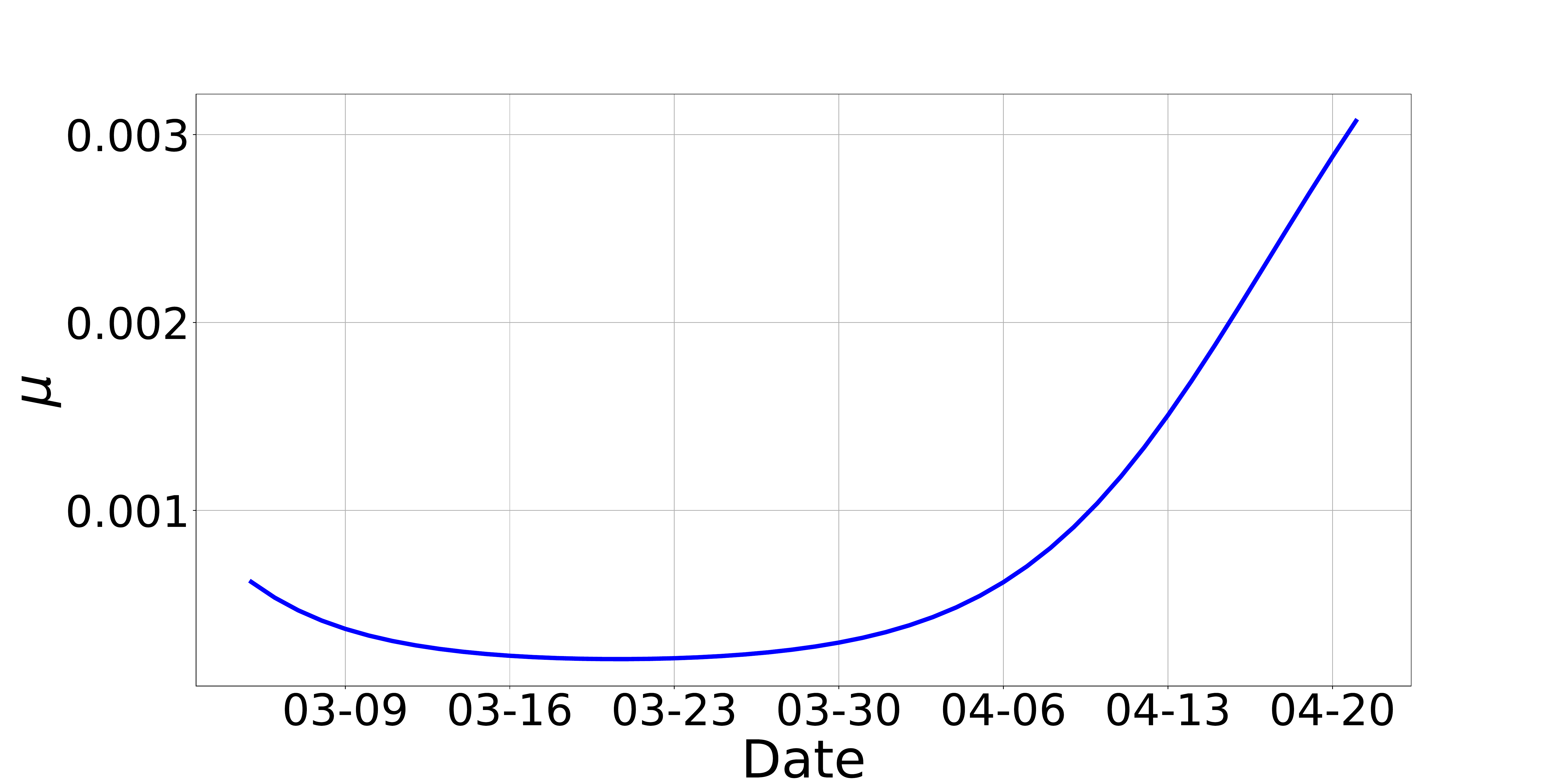}}
\caption{\label{Fig:Inference}fPINNs results: Inferences of unobserved dynamics and model parameters based on the reported data from 27 February 2022 to 21 April 2022.}
\end{figure*}

Next we present inferences of the unobserved dynamics of the compartments $S(t)$ and $E(t)$, the time-independent model parameter $\alpha$, the time-dependent model parameters $\beta(t)$ and $\mu(t)$. The training results are based on the reported data from 27 February 2022 to 21 April 2022. The fractional order in the fractional SEIR model is inferred as $\alpha = 0.75$. Inferences for the time-dependent dynamics are shown in Fig. \ref{Fig:Inference}. 

\subsection{\label{sec:Prediction}Prediction}
As the data fitting and inferences of fPINNs for the fractional SEIR model are already obtained, prediction for the dynamics of COVID-19 caused by Omicron variant can be made. In this part, forecasts of the current infected cases and daily new infected cases are displayed and are compared with the ground truth.  

The fractional SEIR model is a system of nonlinear ODEs. To obtain the prediction results, the initial conditions for all of the compartments as well as the model parameters should be known.  The initial conditions are either from the training data or the inferences from the training time window. The model parameters are calibrated in the training time window. For their values in the prediction time window, the selection strategies for $\alpha$, $\beta(t)$, and $\mu(t)$ are different in this paper. The time-independent parameter $\alpha$ is assumed to be the same as the one in the training time window.  The parameter $\beta(t)$ that indicates the rate of transmission from $S(t)$ compartment to $E(t)$ compartment is a critical parameter determining the trend of the epidemic. In the prediction part, the mean value of $\beta(t)$ is assumed to be its final value $\beta(t_{u}^{N_u})$ of the training time window and an uncertainty bound of 30\% is added to the mean value to show how this uncertainty propagates for the infected cases. For the parameter $\mu(t)$ which denotes the rate of transmission from the $I(t)$ compartment to $R(t)$ compartment, its value is assumed to be its final value $\mu(t_{u}^{N_u})$ of the training time window. Here are two sets of prediction results.  

\begin{figure*}[h] 
\centering
\subfigure[\label{Fig:beta_predict}Model parameter $\beta(t)$ for the one-week prediction based on reported data from 27 Feb 2022 to 21 April 2022]{\includegraphics[width = 8cm]{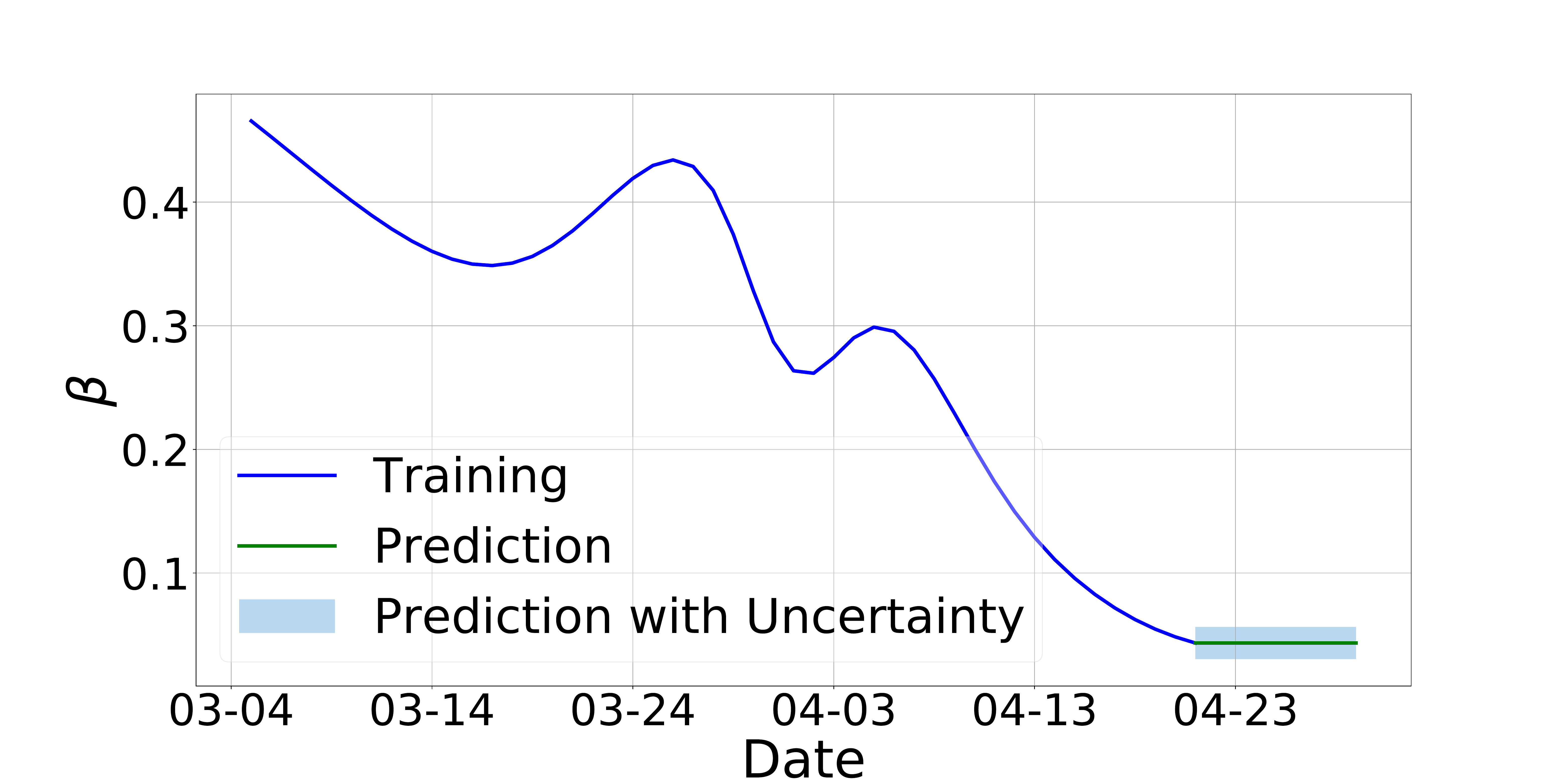}}
\subfigure[\label{Fig:beta_predicts}Model parameter $\beta(t)$ for the one-week prediction based on reported data from 27 Feb 2022 to 30 April 2022]{\includegraphics[width = 8cm]{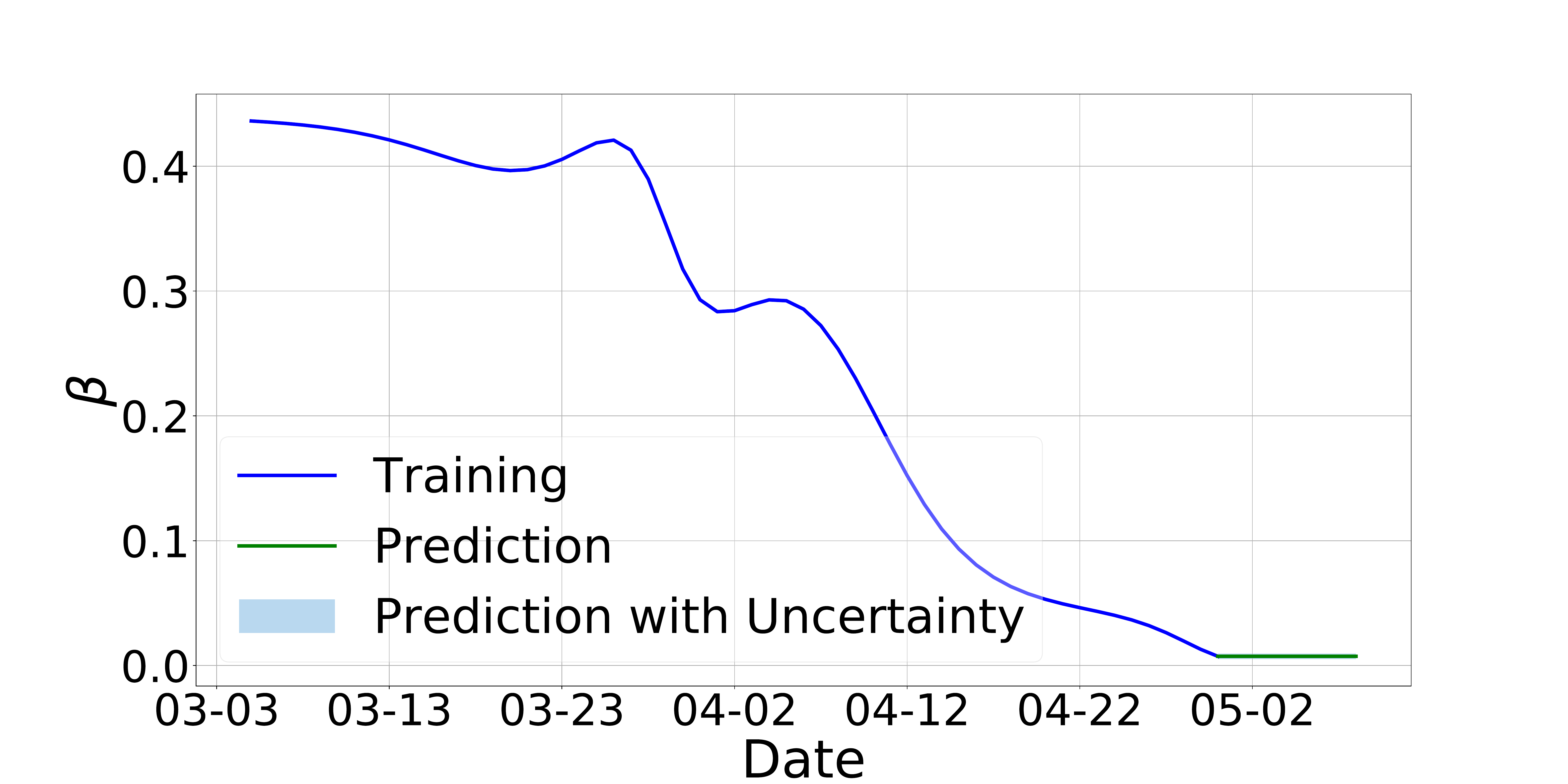}} 
\\
\centering
\subfigure[\label{Fig:Infected_predict}One-week prediction for infected cases based on reported data from 27 Feb 2022 to 21 April 2022: $\alpha = 0.75$, $\mu(t) = 3.07\times10^{-3}$]{\includegraphics[width = 8cm]{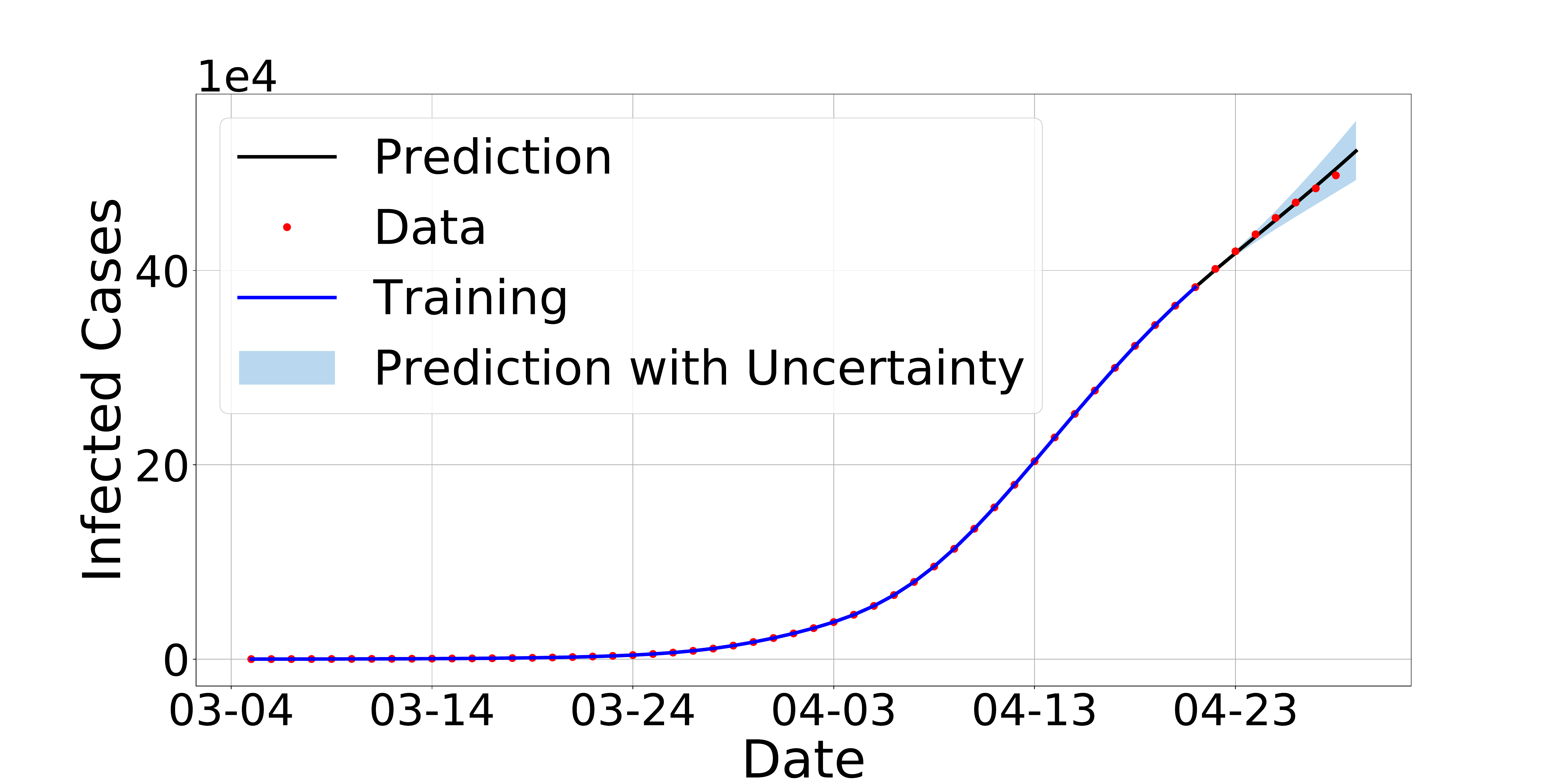}}
\subfigure[\label{Fig:Infected_predicts}One-week prediction for infected cases based on reported data from 27 Feb 2022 to 30 April 2022: $\alpha = 0.96$, $\mu(t) = 5.11\times10^{-3}$]{\includegraphics[width = 8cm]{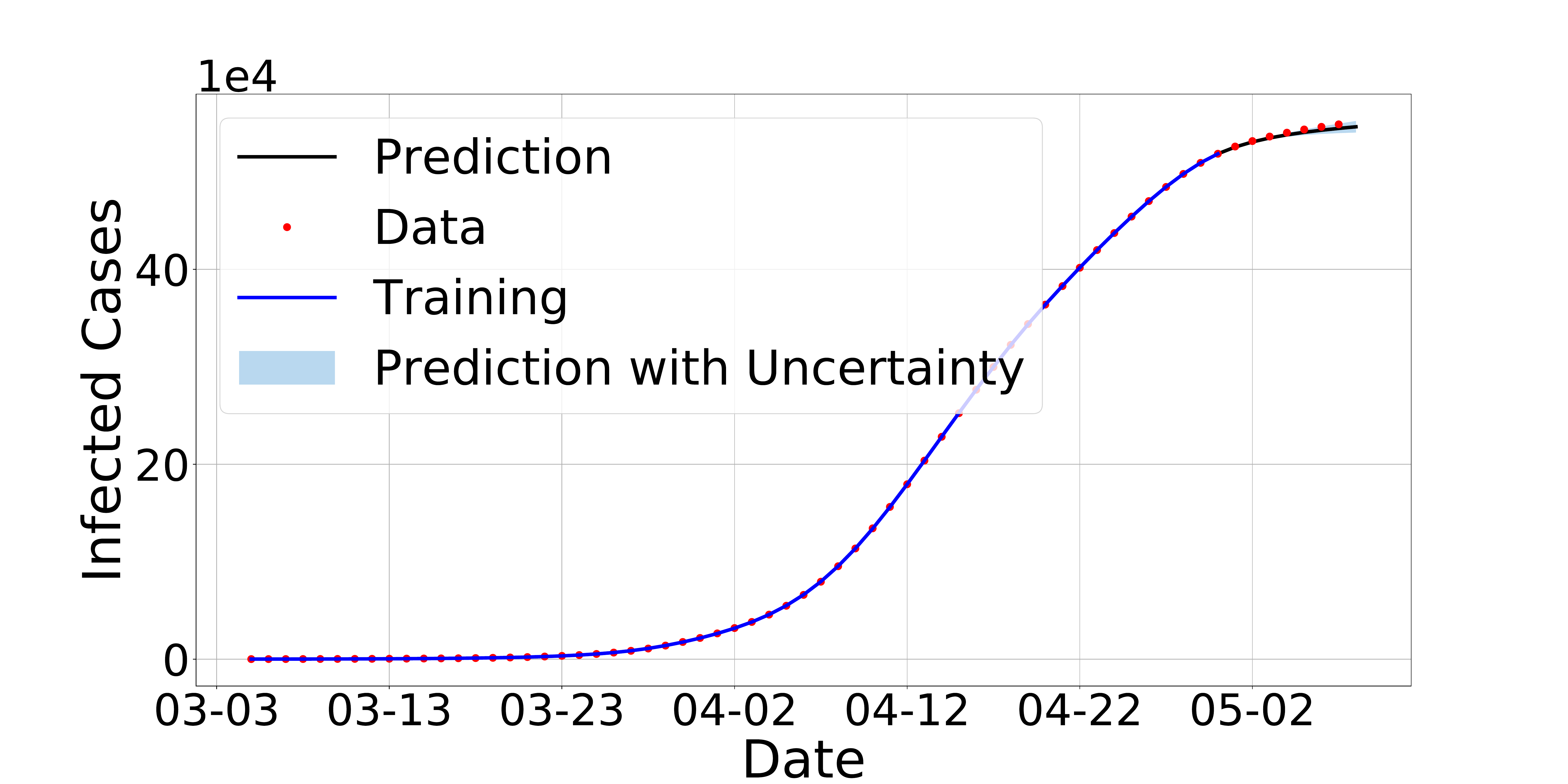}}
\\
\centering
\subfigure[\label{Fig:New_Infected_predict}One-week prediction for new infected cases based on reported data from 27 Feb 2022 to 21 April 2022: $\alpha = 0.75$, $\mu(t) = 3.07\times10^{-3}$]{\includegraphics[width = 8cm]{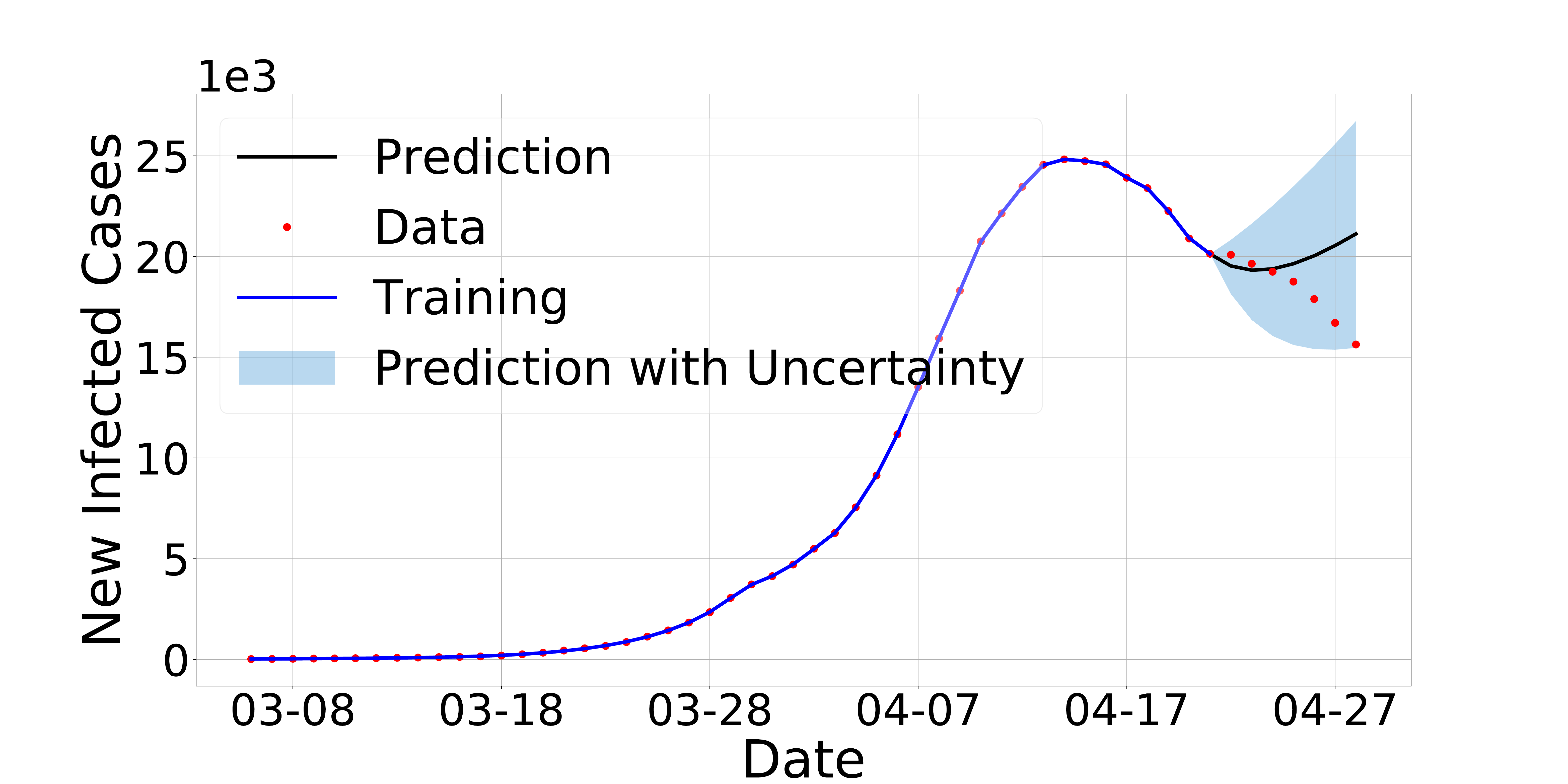}}
\subfigure[\label{Fig:New_Infected_predicts}One-week prediction for new infected cases based on reported data from 27 Feb 2022 to 30 April 2022: $\alpha = 0.96$, $\mu(t) = 5.11\times10^{-3}$]{\includegraphics[width = 8cm]{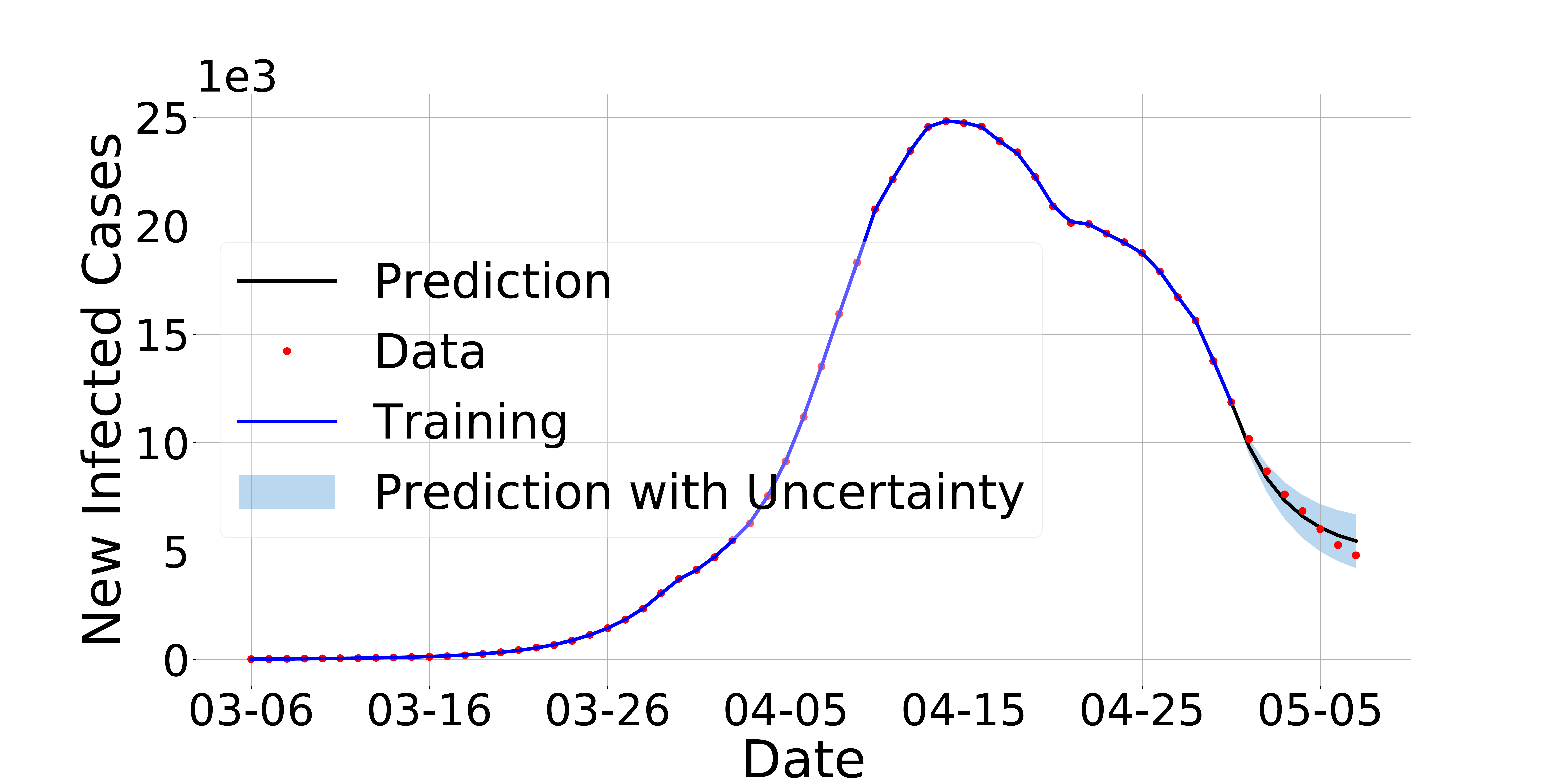}}
\caption{\label{Fig:Prediction}One-week predictions based on the learned fractional SEIR model. Left column: Predictions based on the reported data from 27 Feb 2022 to 30 April 2022; Right column: Predictions based on reported data from 27 Feb 2022 to 30 April 2022. }
\end{figure*}

Fig. \ref{Fig:beta_predict} displays the parameter $\beta(t)$ for the one-week prediction based on reported data from 27 Feb 2022 to 21 April 2022.  Corresponding prediction results on the infected cases and the new infected cases are shown in Fig. \ref{Fig:Infected_predict} and Fig. \ref{Fig:New_Infected_predict}, respectively. It can be observed from Fig. \ref{Fig:New_Infected_predict} that the learned fractional SEIR model in this paper is able to capture the sudden change of the tendency for the new infected cases. In addition, the prediction falls within the uncertainty bounds.

Fig. \ref{Fig:beta_predicts} displays the parameter $\beta(t)$ for the one-week prediction based on reported data from 27 Feb 2022 to 30 April 2022.  Corresponding prediction results on the infected cases and the new infected cases are shown in Fig. \ref{Fig:Infected_predicts} and Fig. \ref{Fig:New_Infected_predicts}, respectively. Data fitting and inferences for this case are given by Fig. \ref{Fig:DataFitting_v2} and Fig. \ref{Fig:Inference_v2} in Appendix \ref{Sec:Appendix}. It is evident that the predictions successfully forecast the dynamics when the ground truth of the new infected cases 
changes monotonously.

\section{\label{sec:Discussion}Conclusion}
In this paper we update the classical SEIR model by introducing the Caputo-Hadamard derivative to describe the relatively slow increase in the detected cases of the daily new infected at the beginning of the COVID-19 outbreak caused by Omicron variant. In order to draw a complete picture of the modified model for the epidemic, we employ fPINNs to calibrate the unknown model parameters and unobserved dynamics. Furthermore, reliable predictions are made to show the feasibility and predictability of the fPINNs for the fractional SEIR model. The framework of modeling and predicting the transmission of COVID-19 caused by Omicron variant is also applicable for other areas around the world, such as Berlin and New York.  

\begin{acknowledgments}
This work is supported by the National Natural Science Foundation of China (No. 11872234).
\end{acknowledgments}

\section*{Author Declarations}
\subsection*{Conflict of Interest}
The authors declare that they have no conflicts of interest to disclose.
\subsection*{Author Contributions}
The authors equally contributed to this work.

\section*{Data Availability Statement}
The original data of the simulations in the present paper are collected from "Epidemic Notification (in Chinese)" published by National Health Commission of the People’s Republic of China \cite{EpidemicNotification}. 
The data sets generated during the current study are available from the corresponding author upon reasonable request.

\appendix
\section{\label{Sec:Appendix} Data Fitting and Inferences based on Reported Data from 27 Feb 2022 to 30 April 2022}
\begin{figure*}[h]
\centering
\subfigure[\label{Fig:New_I_v2}New infected cases $I^{n}(t)$]{\includegraphics[width = 8cm]{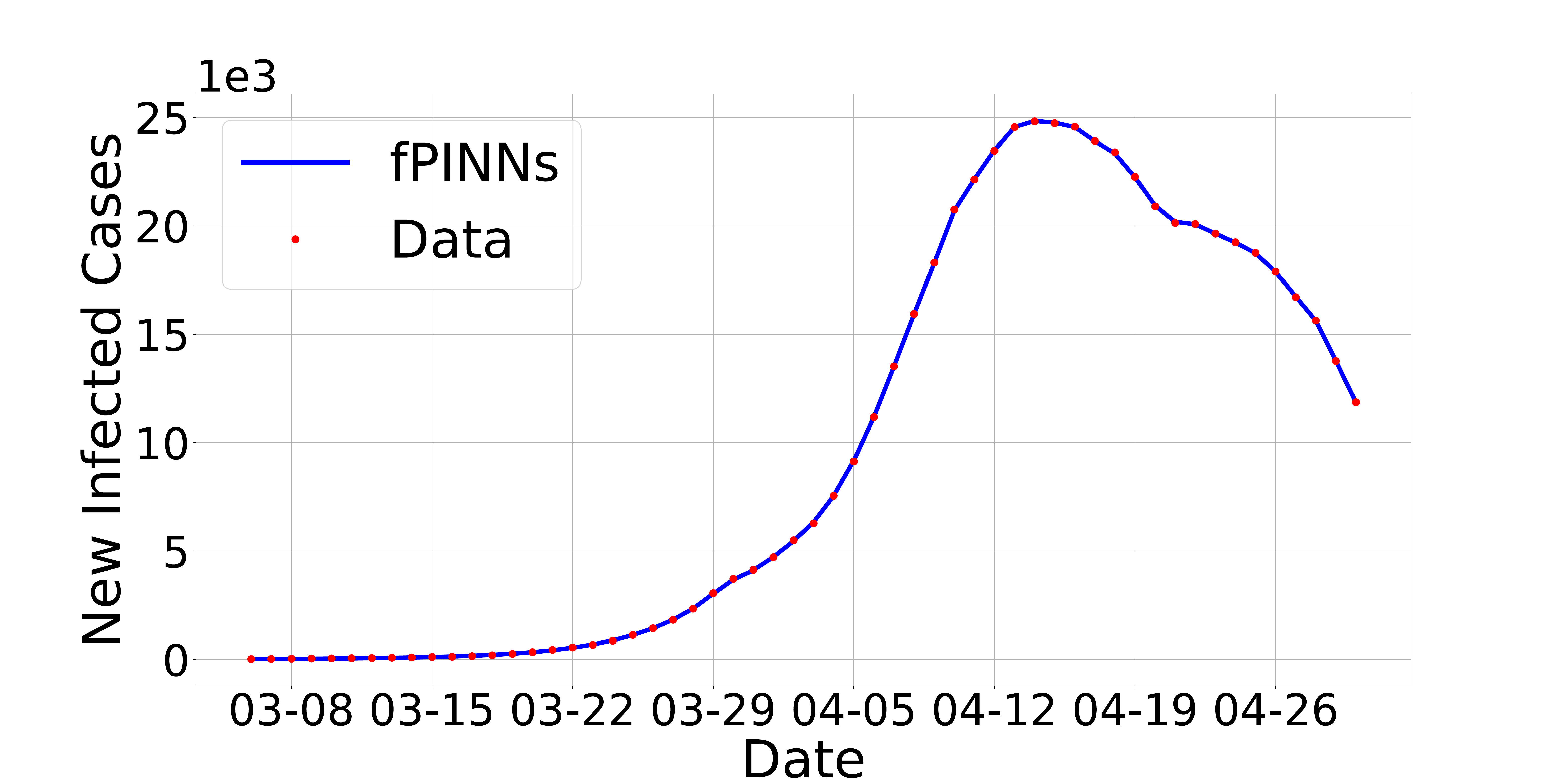}}
\subfigure[Cumulative infected cases $I^{c}(t)$]{\label{Fig:Cumulative_I_v2}\includegraphics[width = 8cm]{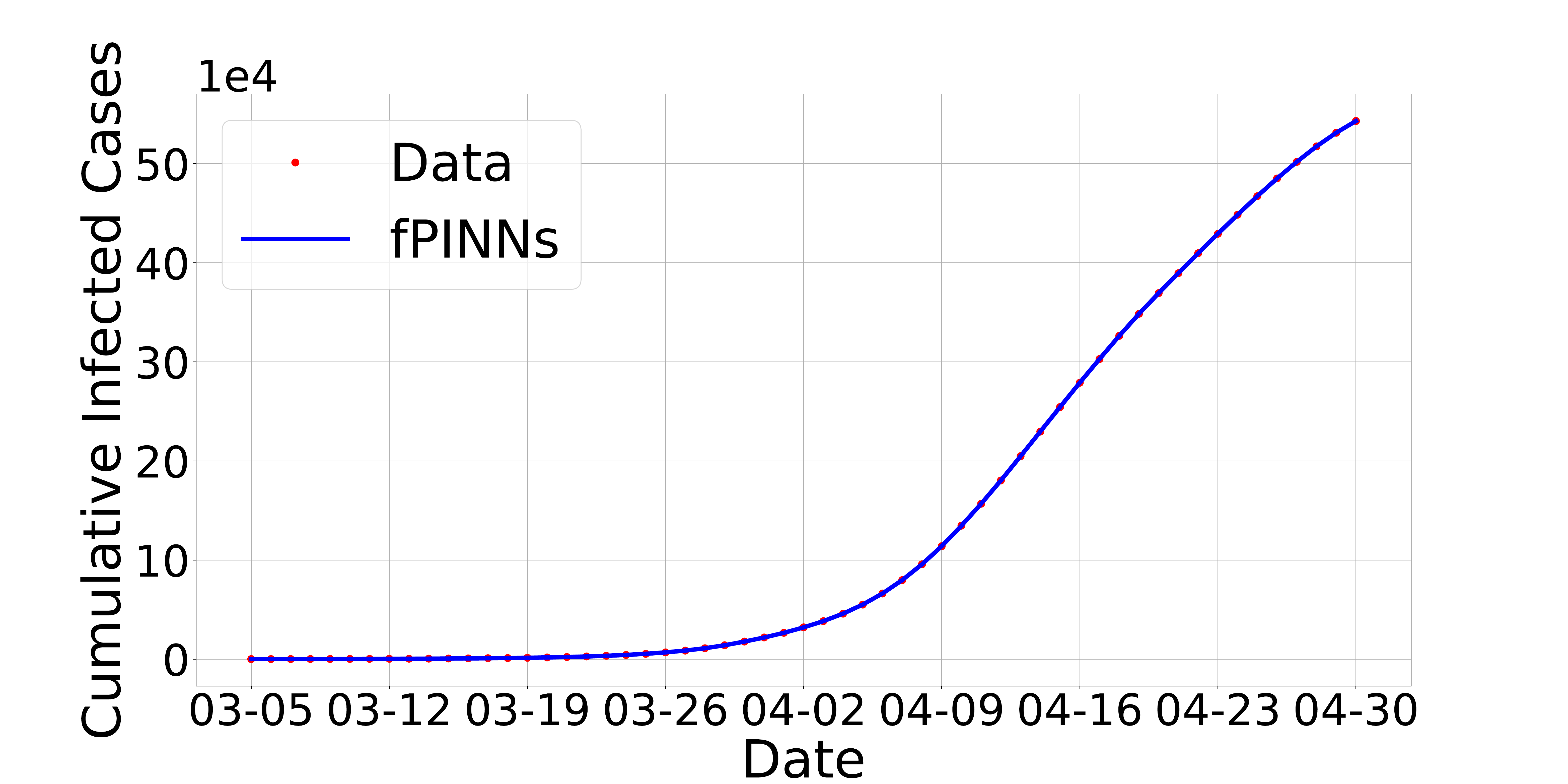}}
\\
\centering
\subfigure[\label{Fig:New_R_v2}New removed cases $R^{n}(t)$]{\includegraphics[width = 8cm]{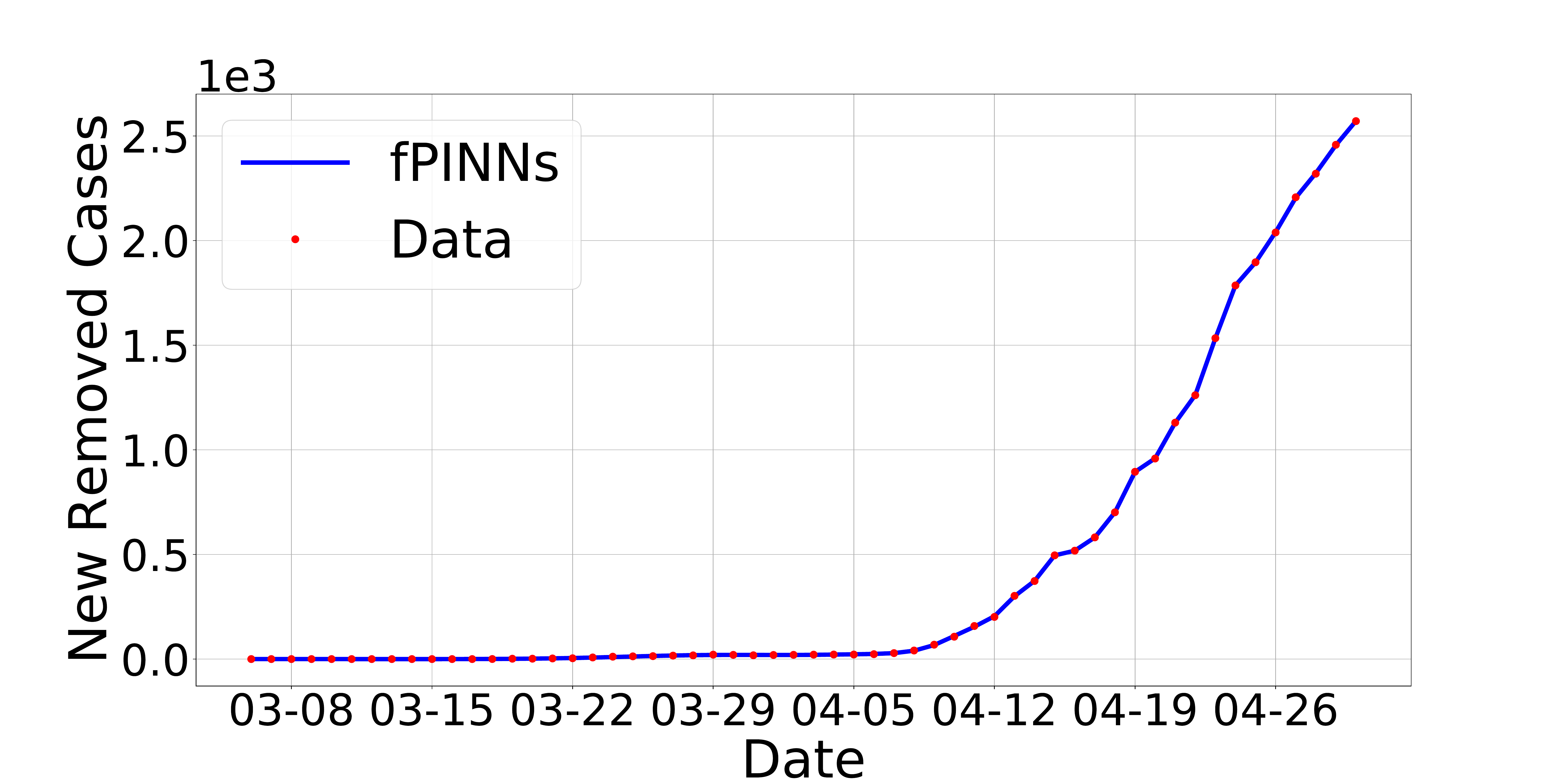}}
\subfigure[\label{Fig:Current_R_v2}Current removed cases $R(t)$]{\includegraphics[width = 8cm]{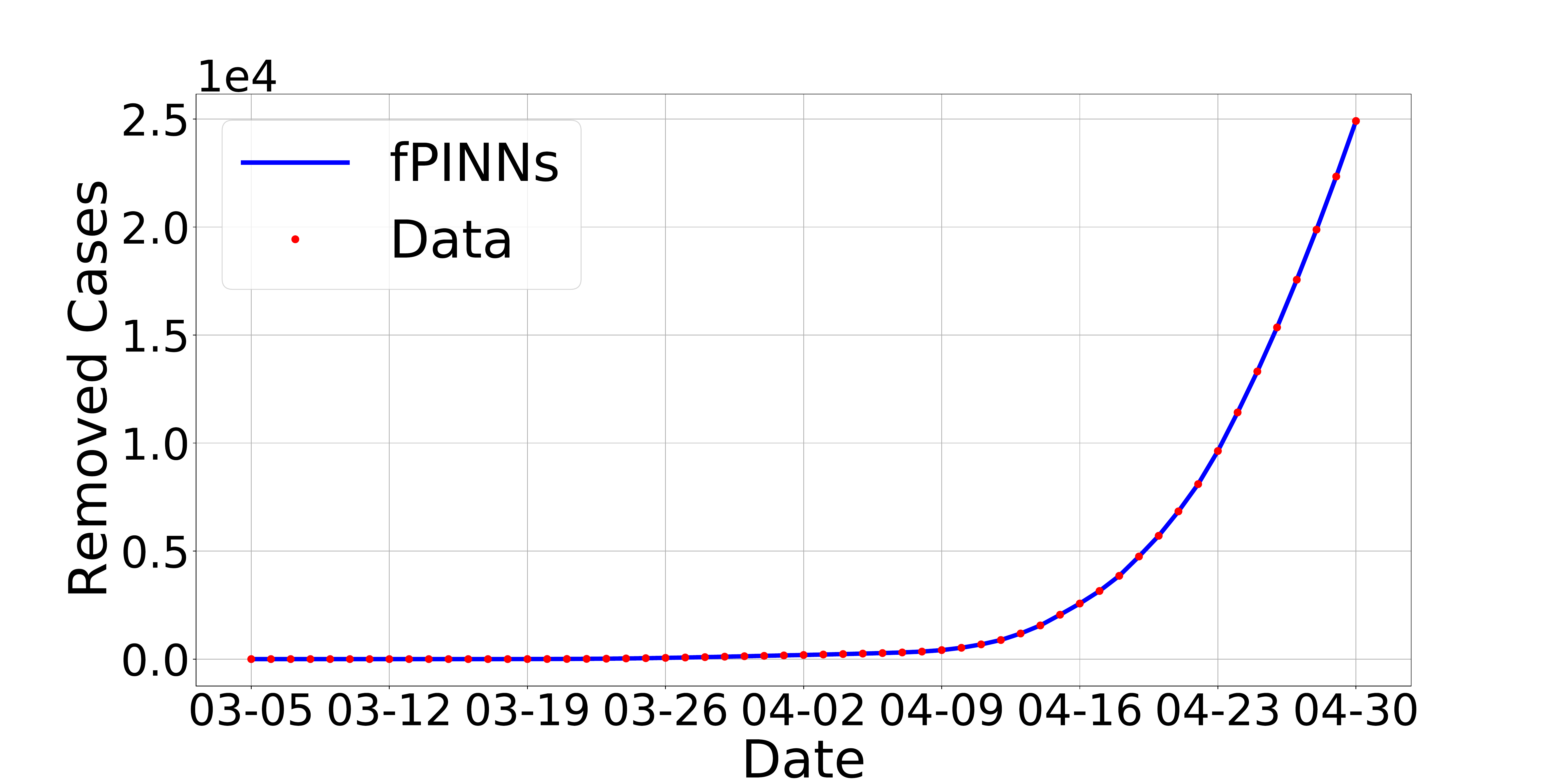}}
\\
\centering
\subfigure[\label{Fig:Current_I_v2}Current infected cases $I(t)$]{\includegraphics[width = 8cm]{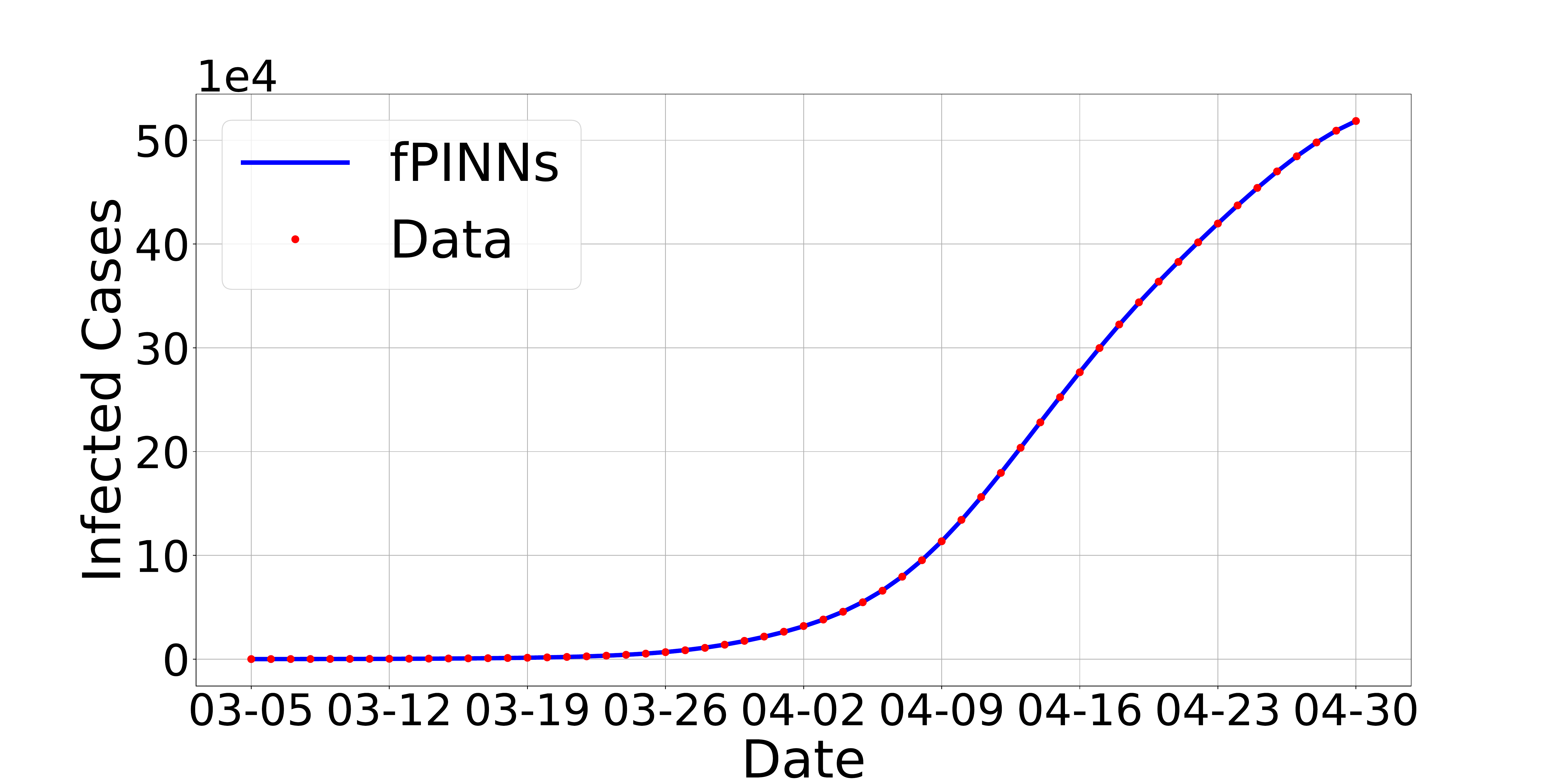}}
\caption{\label{Fig:DataFitting_v2}fPINNs results: Data fitting based on the reported data from 27 February 2022 to 30 April 2022.}
\end{figure*}

\begin{figure*}[h]
\centering
\subfigure[\label{Fig:Current_S_v2}Current susceptible cases $S(t)$]{\includegraphics[width = 8cm]{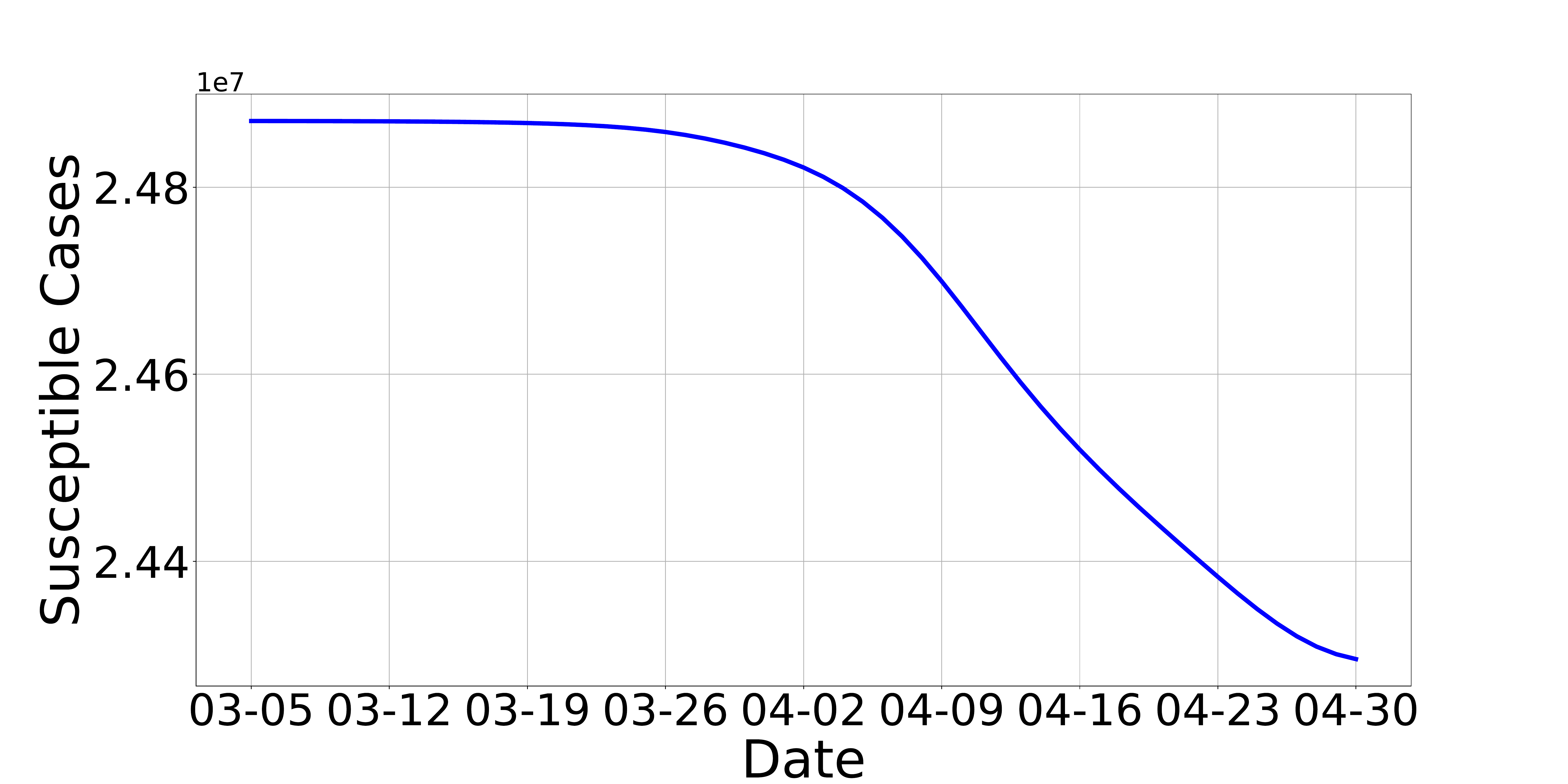}}
\subfigure[\label{Fig:Current_E_v2}Current exposed cases $E(t)$]{\includegraphics[width = 8cm]{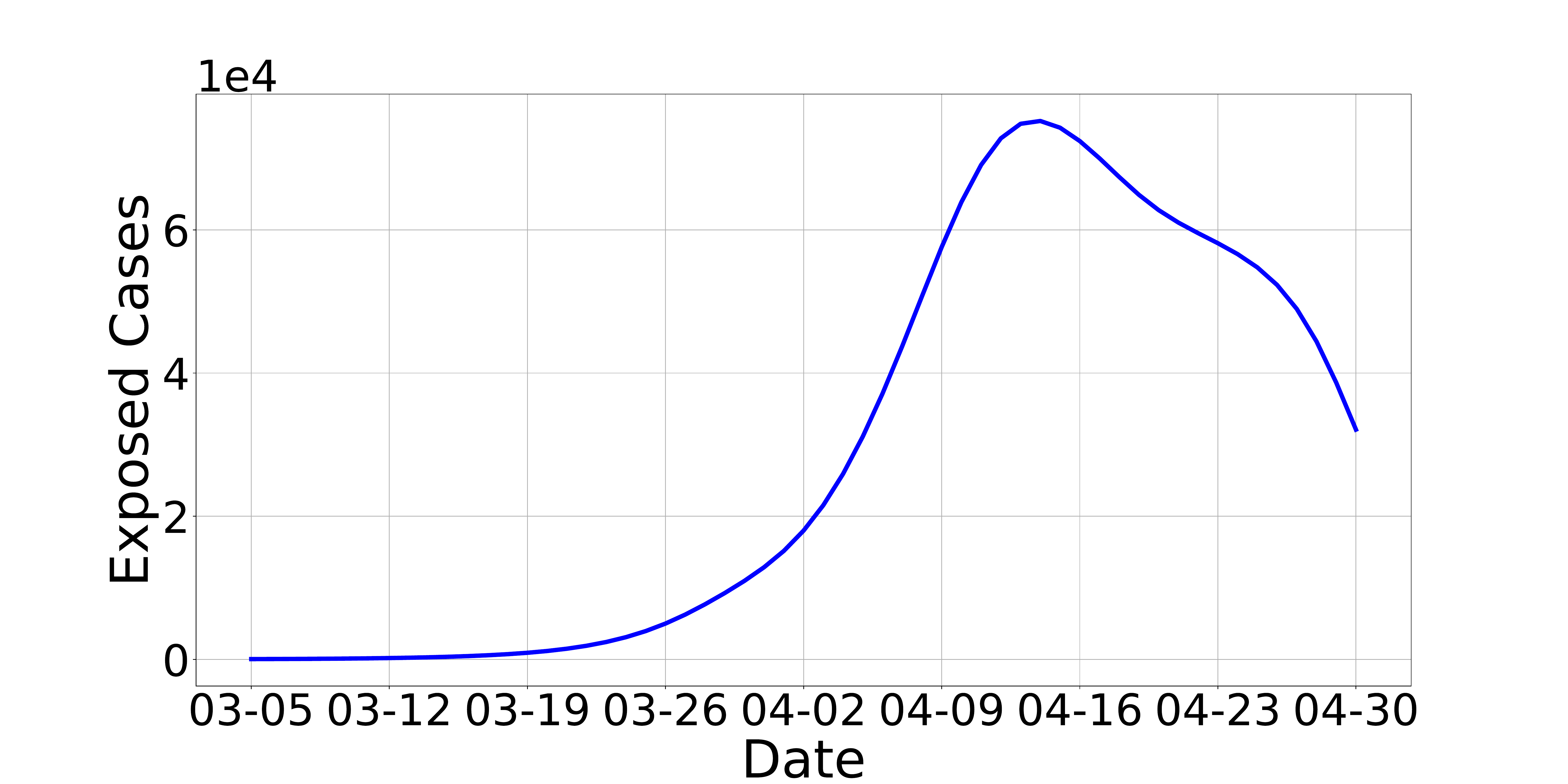}}
\\ 
\centering
\subfigure[\label{Fig:Beta_v2}Model parameter $\beta(t)$]{\includegraphics[width = 8cm]{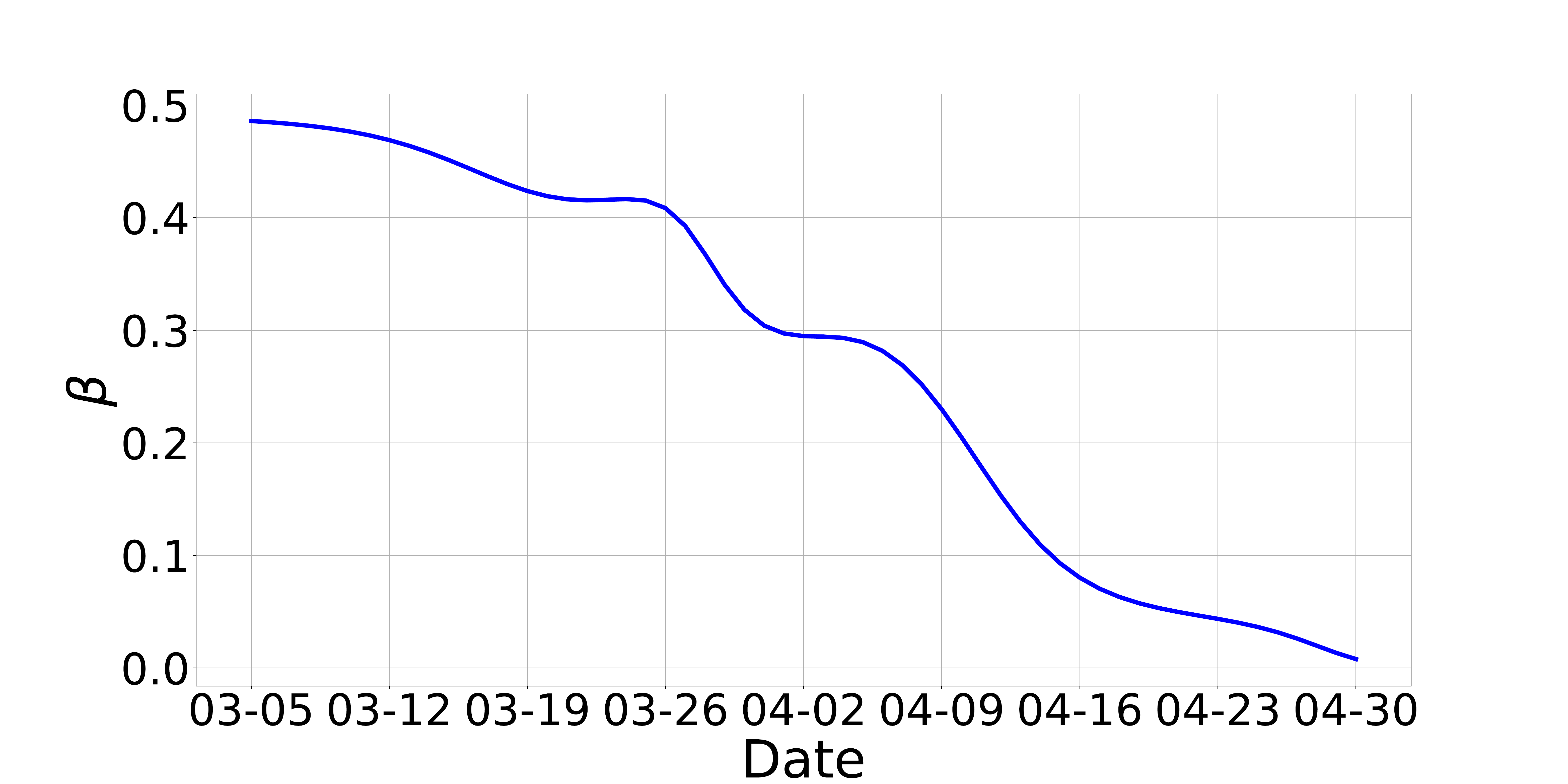}}
\subfigure[\label{Fig:mu_v2}Model parameter $\mu(t)$]{\includegraphics[width = 8cm]{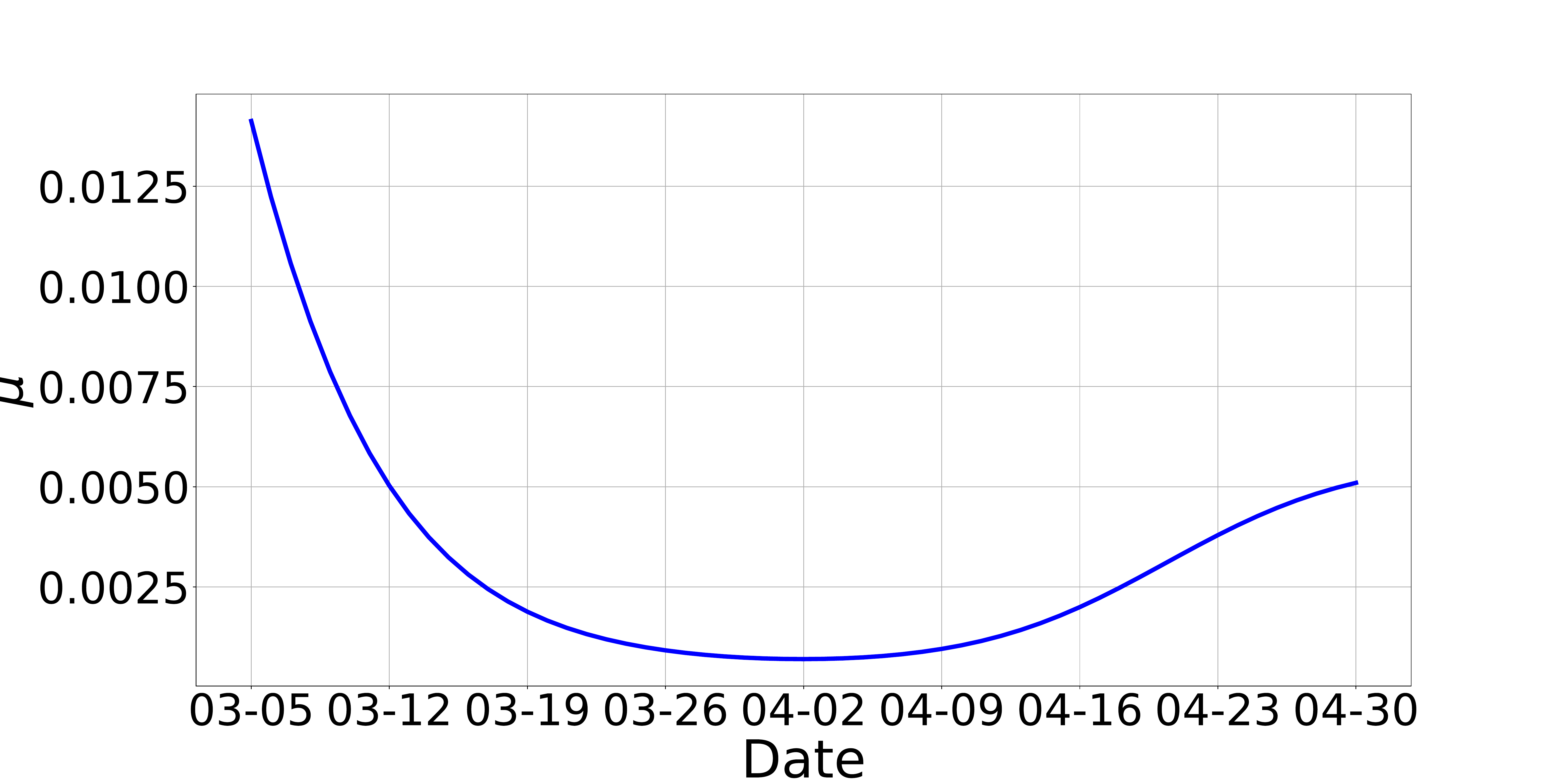}}
\caption{\label{Fig:Inference_v2}fPINNs results: Inferences of unobserved dynamics and model parameters based on the reported data from 27 February 2022 to 30 April 2022.}
\end{figure*}

In this section, data fitting and inferences based on reported data of Shanghai from 27 Feb 2022 to 30 April 2022 are shown in Fig. \ref{Fig:DataFitting_v2} and Fig. \ref{Fig:Inference_v2}. The corresponding training data is therefore from 5 March 2022 to 30 April 2022 after taking seven-days average of the reported data. In this case, the fractional order is inferred as $\alpha = 0.96$.

\section*{References}

\end{document}